\documentclass[11pt]{article}

\usepackage{amsmath}
\usepackage{amssymb}
\usepackage[all]{xy}
\usepackage{mathrsfs}

\DeclareMathAlphabet{\scr}{U}{eus}{m}{n}

\newcommand\N{{\mathbb N}}
\newcommand\Z{{\mathbb Z}}
\newcommand\Q{{\mathbb Q}}
\newcommand\F{{\mathbb F}}

\newcommand\A{{\mathbb A}}
\newcommand\G{{\mathbb G}}

\newcommand\cA{{\cal A}}
\newcommand\cB{{\cal B}}

\newcommand\cF{{\cal F}}

\newcommand\cH{{\cal H}}

\newcommand\cN{{\cal N}}

\newcommand\cQ{{\cal Q}}

\newcommand\xB{B}

\newcommand\cNN{{\cal N}^-}

\newcommand\sH{{\mathscr H}}

\newcommand\Mloc{{\rm\bf M}^{\rm loc}}

\newcommand\Adm{{\rm Adm}}
\newcommand\Tr{{\rm Tr}}
\newcommand\Frob{{\rm Fr}}
\newcommand\Fr{{\rm Fr}}
\newcommand\xt{[\hspace{-.12em}[ t ]\hspace{-.12em} ] } 
\newcommand\xT{(\!( t )\!)}

\newcommand\Ql{\overline{\mathbb Q}_\ell}

\newcommand\Fl{\mathcal Fl}
\newcommand\Wext{\widetilde{W}}

\newcommand\tens\otimes

\newcommand\lto{\longrightarrow}
\newcommand\id{\mathop{\rm id}\nolimits}

\newcommand{\Spec}{\mathop{\rm Spec}}

\newcommand{\Grass}{\mathop{\rm Grass}\nolimits}

\newcommand\Gal{\mathop{\rm Gal}}

\newcommand\gfrac[2]{\genfrac{}{}{0pt}{}{#1}{#2}}

\newcommand\qed{\hfill$\square$}

\newtheorem{thm}{Theorem}[section]
\newtheorem{stz}[thm]{Proposition}
\newtheorem{lem}[thm]{Lemma}
\newtheorem{Def}[thm]{Definition}
\newtheorem{kor}[thm]{Corollary}

\author{Ulrich G\"{o}rtz and Thomas J. Haines \footnote{Research supported in part by NSF grant DMS 0303605}}
\title{The Jordan-H\"older series for nearby cycles \\
on some Shimura varieties and affine flag varieties}
\date{}

\begin{document}

\maketitle

\bigskip

\section{Introduction}

This article has several aims.  The initial motivation was our desire to understand, as explicitly as possible, the irreducible constituents of the nearby cycles sheaf $R\Psi$ on the reduction modulo $p$ of a Shimura variety with Iwahori-level structure at a fixed prime $p$.  In some sense, the complexity of $R\Psi$ is a measure of the singularities in the reduction modulo $p$.  Moreover, we were motivated by the role $R\Psi$ plays in the computation of the semi-simple local zeta function at $p$ for such a Shimura variety, cf. \cite{RZ1}, \cite{HN1}.  

Via the relationship between the Shimura variety and its Rapoport-Zink local 
model \cite{RZ2}, one can translate the problem into that of understanding the nearby cycles on a local model.  The latter can be embedded as a finite-dimensional piece of a $\Z_p$-ind-scheme $M$ which is a deformation of the affine Grassmannian ${\rm Grass}_{{\mathbb Q}_p}$ to the affine flag variety $\Fl_{\F_p}$ for the underlying $p$-adic group $G$ ($M$ exists at least if $G$ is either $GL_n$ or $GSp_{2n}$; see \cite{HN1} and section 8).  A very similar deformation ${\rm Fl}_X$ over a smooth curve $X$ (due to Beilinson) exists for any group $G$ in the function field setting, and has been extensively studied by Gaitsgory \cite{Gaitsgory}.  

The ``maximal parahoric'' subgroup $G(\mathcal O)$ in the loop group of $G$ acts on ${\rm Grass}$, the generic fiber of the deformation $M$ or ${\rm Fl}_X$.  Fixing a dominant coweight $\mu$ of $G$, one may therefore consider nearby cycles $R\Psi = R\Psi(IC_\mu)$, 
taken with respect to the deformation, where $IC_\mu$ is the intersection complex on the closure $\overline{\cQ}_\mu$ of the $G(\mathcal O)$-orbit 
in ${\rm Grass}$ indexed by $\mu$.  Letting $\cB$ denote the standard Iwahori subgroup of $G(\mathcal O)$, one may show that $R\Psi$ is a $\cB$-equivariant perverse sheaf on the affine flag variety, with certain additional properties (e.g., it belongs to the category $P^{\cB}_q(\Fl,\Ql)$ of section 4).  Thus, quite generally, we were led to the problem of understanding the Jordan-H\"{o}lder series for
objects $\cF$ of the category $P^{\cB}_q(\Fl,\Ql)$.  It turns out that the only irreducible objects in this category are Tate-twists of intersection complexes on finite-dimensional Schubert varieties $\cB_{\overline{w}}$ in $\Fl$ (they are indexed by elements $w$ in the extended affine Weyl group $\widetilde{W}$ for $G$; we denote such intersection complexes by $IC_w$).  Thus, we may define non-negative ``multiplicities'' $m(\cF,w,i)$ by the identity in $P^{\cB}_q(\Fl,\Ql)$
$$
\cF^{ss} = \bigoplus_{w \in \widetilde{W}} \bigoplus_{i \in \Z} IC_w (-i)^{\oplus m(\cF, w, i)}.
$$

Our main theorem is the following cohomological interpretation for the
integers $m(\cF,w,i)$ for sheaves $\cF$ which satisfy a somewhat technical
hypothesis (Property (P) of section 4.3). This property holds for all the nearby cycles we consider, and also for a much larger class of objects (those with a suitable filtration by Wakimoto sheaves; see section 7).  

\begin{thm} Suppose $\cF \in P^{\cB}_q(\Fl,\Ql)$ satisfies property (P) for the integer $d$.  Then for each $w \in \widetilde{W}$,
$$
\sum_i m(\cF, w,i)q^i =  (-q)^d \Tr(\Fr_q, H_c^\bullet(\Fl, D\cF \tens
IC(\cB^{\overline{w}}))).
$$
\end{thm}
Here $D\cF$ denotes the Verdier dual of $\cF$; there is also a ``dual''
formulation (see Theorem \ref{cohomological interpretation}) in which $\cF$
appears rather than $D\cF$.  We introduce and study the intersection complex
$IC(\cB^{\overline{w}})$ appearing here in section 3.  In a certain sense, it
is the intersection complex  for an orbit $\cB^w$ for the ``opposite'' Iwahori
subgroup $\cB^-$.  Such orbits are infinite-dimensional, and so it is
necessary to construct $IC(\cB^{\overline{w}})$ by pull-back from certain
finite-dimensional quotients (introduced by Faltings \cite{F}) of open sets
$\Omega \subset \Fl$.  The basic properties of $IC(\cB^{\overline{w}})$ play a key role in the proof.  The other main ingredient is the cohomological interpretation of inverse Kazhdan-Lusztig polynomials, announced in \cite{KL2}.  In fact, the Kazhdan-Lusztig theorem can be seen as a special case of our main theorem (see section 4.5).  Since no proof of the Kazhdan-Lusztig theorem in exactly our setting has appeared in the literature, we provide a proof for this result in section 3.  Here again we make essential use of our sheaf $IC(\cB^{\overline{w}})$ and of the quotients of Faltings \cite{F} to pass from ind-schemes to ordinary finite-dimensional schemes.

The main theorem was discovered in the attempt to find a conceptual explanation for the following remarkable fact (Cor. \ref{remarkable_fact}): Let $M_\mu$ denote the scheme-theoretic closure in our deformation of $\overline{\cQ}_\mu$.  Let $\tau$ denote the element of $\widetilde{W}$ indexing the unique dimension zero Iwahori-orbit in the 
special fiber of $M_\mu$ (i.e., the ``most singular point'').  Then we have

\begin{kor} For $R\Psi = R\Psi(IC_\mu)$, the multiplicity $m(R\Psi, \tau, i)$ is the $2i$-th intersection Betti number for $\overline{\cQ}_\mu$.  
\end{kor}
This fact, noticed after several computer-aided calculations using the algorithm described below, gave us the first hint that multiplicity functions have a cohomological interpretation.

It is clear that knowing the multiplicities $m(\cF,w,i)$ for every $w$ and every $i$ is equivalent to knowing an explicit formula for the trace function $\Tr(\Fr_q,\cF)$ as an element in the Iwahori-Hecke algebra (see section 4.2).  On the other hand, the Kottwitz conjecture (proved in \cite{Gaitsgory} and \cite{HN1}, cf. section 2.7) identifies the semi-simple trace function $\Tr^{ss}(\Fr_q,R\Psi(IC_\mu))$ as an explicit sum of functions $\Theta_\lambda$, each of which can be computed (with a computer program, see sections 7,9).  Thus, there is an algorithm to compute the multiplicities $m(R\Psi, w, i)$ in any given case.  (Actually, to justify the algorithm, one has to know that $\Tr^{ss}(\Fr_q,R\Psi) = \Tr(\Fr_q,R\Psi)$; see comments below.)  Similiarly, the algorithm works to compute the numbers $m(\cF,w,i)$ whenever $\cF$ has a suitable filtration by Wakimoto sheaves, cf. section 7.  In section 9, we 
give the results of several such computer-aided calculations, and we give explanations for some empirical observations we made.

Another consequence of the Kottwitz conjecture is that the functions $\Tr^{ss}(\Fr_q, R\Psi(IC_\mu))$ form a basis for the center $Z(\mathcal H)$ of the Iwahori-Hecke algebra for $G$, as $\mu$ ranges over dominant coweights.  In terms of the standard generators $T_x = {\rm char}(\cB x \cB)$ ($x \in \widetilde{W}$) for ${\mathcal H}$, we have the identity
$$  \Tr^{ss}(\Fr_q, R\Psi) = \sum_x \left( \sum_{w \geq x} \left( \sum_i m(R\Psi,w,i)q^i \, \varepsilon_w P_{x,w} \right) \right) \, T_x, $$
where $\leq$ is the Bruhat order on $\widetilde{W}$ and $P_{x,w}$ are the Kazhdan-Lusztig polynomials \cite{KL1}, cf. section 4.2.  
Our results on $m(R\Psi,w,i)$ therefore give us both a conceptual and an algorithmic way to make this formula explicit, and this works for any group $G$.
Other approaches to computing elements in $Z(\mathcal H)$ explicitly in terms of the $T_x$-basis were given previously in \cite{H1},\cite{H3}, \cite{HP}, and \cite{S}.  

Other aims of the paper were to provide proofs of several necessary results of a foundational nature, some of which are very general and might be useful in other contexts.  Here is a list of the highlights:

\noindent $\bullet$  A study of categories of $\cB$-equivariant perverse sheaves on $\Fl$; section 4;

\noindent $\bullet$  For finite-type schemes $X$ over a trait $(S,s,\eta)$, the canonical decomposition of the category $P(X \times_s \eta,\Ql)$ of middle perverse sheaves on $X_{\bar{s}}$ endowed with a compatible continuous action of 
${\rm Gal}(\bar{\eta}/\eta)$, into ``unipotent'' and ``non-unipotent'' subcategories; section 5;

\noindent $\bullet$  The proof of Gaitsgory's theorem that our sheaves $R\Psi$ are unipotent (in the above sense), and the consequence that $\Tr^{ss}(\Fr_q,R\Psi) = \Tr(\Fr_q, R\Psi)$, alluded to above; section 5;

\noindent $\bullet$  If $X$ is finite-type over $S$, an unequal characteristic Henselian trait, and $X_\eta$ is smooth, then the nearby cycles sheaf $R\Psi^X(\Ql)$ is {\em mixed}; section 10;

\noindent $\bullet$  If $X$ as above is proper, the existence of the weight spectral sequence; section 8.

Returning now to the situation of Shimura varieties $Sh$ considered in section 8, the Shimura datum $(G,X)$ provides us with a minuscule coweight $\mu$ of $G_{{\mathbb Q}_p}$, and the model $M_\mu$ above coincides with the Rapoport-Zink local model ${\bf M}^{\rm loc}$, so the above results apply to give the Jordan-H\"{o}lder series for $R\Psi_{Sh}$.  
The irreducible constituents are Tate-twists of intersection complexes $IC_{Sh,w}$ relative to a stratification $\coprod Sh_{w}$ on the special fiber $Sh_{\F_p}$, induced by the $\cB$-orbit stratification on $M_\mu$ (see \cite{GN}), indexed by the subset ${\rm Adm}(\mu)$ of section 2.4 
\footnote{Recall that $Sh$ is usually geometrically disconnected.  However the various connected components carry the ``same'' stratification, induced in each case by that on $M_\mu$.  To be precise, in all statements the symbol $IC_{Sh,w}$ should be interpreted as the direct sum, over all geometric connected components, of the intersection complexes indexed by $w$ on the various components.  In particular this applies to the results in section 8.}.  
If $Sh$ is proper over the ring of integers $\mathcal O_E$ of the reflex field $E$, then the weight spectral sequence $_WE^{\bullet, \bullet} \Rightarrow H^\bullet(Sh_{\overline{\mathbb Q}_p}, \Ql[\ell(\mu)])$ exists (by virtue of the Appendix) and its $E_1$-term  can be made somewhat explicit.  We have the following result (see section 8 for details).

\begin{thm}
\begin{enumerate}
\item[(a)] In the category $P(Sh \times_s \eta, \Ql)$,  we have 
$$R\Psi^{ss}_{Sh}(\Ql[\ell(\mu)]) = \bigoplus_{w \in {\rm Adm}(\mu)} \bigoplus_{i=0}^{\ell(\mu)-\ell(w)} IC_{Sh,w}(-i)^{\oplus m(R\Psi_{M_\mu},w,i)}.$$
\item[(b)] Assume $Sh$ is proper over $\mathcal O_E$.  Then there is an
  isomorphism of $\Ql$-spaces
$$_WE_1^{pq} = \bigoplus_{\gfrac{w \in {\rm Adm}(\mu)}{\ell(w) + 2j = -p}} 
IH^{p+q+\ell(w)}(\overline{Sh}_w, \Ql)^{\oplus m(R\Psi_{M_\mu},w,j)}.$$
\end{enumerate}
\end{thm}

\bigskip

The multiplicity functions $\sum_i m(R\Psi(IC_\mu),w,i) q^i$ appear to carry fundamental information.  In all the cases we computed, they are actually polynomials in $q$ of degree exactly $\ell(\mu) - \ell(w)$ (and they vanish if $w \notin {\rm Adm}(\mu)$).  It is perhaps striking that we do not know how to prove the polynomial nature of the multiplicity function in general.  It can be proved when $\mu$ is either minuscule or a coweight for $GL_n$, by using the existence of minimal expressions for the functions $\Theta_\lambda$ that arise in the Kottwitz conjecture, see section 9.  We provide a possible geometrical approach in Part II of the Appendix (section 10). 

\section{Notation and preliminaries}

\subsection{Galois structures on derived categories}

As usual, $p$ will denote a fixed prime number, and $q$ will always
denote a power of $p$.  We will work over Henselian traits $(S,s, \eta)$ of residue characteristic $p$.  These will usually be of the form $S= {\rm Spec}(\F_p[[t]])$ (the function field setting), or 
$S= {\rm Spec}(\Z_p)$ (the $p$-adic setting).  We choose a separable closure $\bar{\eta}$ of $\eta$ and define the Galois group $\Gamma = {\rm Gal}(\bar{\eta}/\eta)$ and the inertia subgroup $\Gamma_0 = 
{\rm ker}[{\rm Gal}(\bar{\eta}/\eta) \rightarrow {\rm Gal}(\bar{s}/s)]$, where $\bar{s}$ is the residue field of the normalization $\bar{S}$ of $S$ in 
$\bar{\eta}$.

Suppose the residue field $k(s)$ has cardinality $q$.  Any element of $\Gamma$ which lifts the inverse ${\rm Frob}_q$ of the topological generator $a \mapsto a^q$ of ${\rm Gal}(\bar{s}/s)$ will be called a {\em geometric Frobenius} and will be denoted by $\Fr_q$.  We will often use the same symbol $\Fr_q$ to denote the Frobenius automorphism ${\rm id}_{X_0} \times {\rm Frob}_q: X_0 \times {\rm Spec}(\overline{\F}_q) \rightarrow X_0 \times {\rm Spec}(\overline{\F}_q)$ for an $\F_q$-scheme $X_0$.  For the most part, we leave it to context to dictate what is meant by the symbol $\Fr_q$ in each situation.

Let $X$ be a scheme of finite type over a finite (or algebraically
closed) field $k$. (The following also works if we assume that $k$ is the fraction field of a discrete valuation ring $R$ with finite residue field, and that $X$ is finite-type over $R$, cf. \cite{Ma}.)  Denote by $\overline{k}$ an algebraic closure of $k$, and by $X_{\overline{k}}$ the base change.

We denote by $D^b_c(X, \Ql)$ the 'derived' category of $\Ql$-sheaves on $X$.
Note that this is not actually the derived category of the category of
$\Ql$-sheaves, but is defined via a limit process. See \cite{BBD} 2.2.14
or \cite{Weil2} 1.1.2 for more details; see also section \ref{inertia_action}. 
Nevertheless, $D^b_c(X, \Ql)$ is a triangulated category which admits the
usual functorial formalism, and which can be equipped with a 'natural'
 $t$-structure having as its core the category of $\Ql$-sheaves.
If $f: X \lto Y$ is a morphism of schemes of finite type over $k$, we
denote by $f_\ast$, $f_!$ the derived functors
 $D^b_c(X, \Ql) \lto D^b_c(Y, \Ql)$.

We will denote by $P(X, \Ql)$ the full subcategory of $D^b_c(X, \Ql)$ consisting of middle perverse sheaves.  See \cite{BBD} or \cite{KW} for a detailed discussion of this notion.

Now let $X$ be a scheme of finite type over $s$. The absolute Galois group
 $\Gal(\overline{\eta}/\eta)$ acts on $X_{\bar{s}}$ through its quotient
 $\Gal(\overline{s}/s)$.

\begin{Def} {\rm (cf. \cite{SGA 7}, exp. XIII)}
\begin{enumerate}
\item[(1)]
The category $D^b_c(X \times_{s} \eta, \Ql)$ is the category of sheaves
$\cF \in D^b_c(X_{\bar{s}}, \Ql)$ together with a continuous action of 
 $\Gal(\overline{\eta}/\eta)$ which is compatible with the action on $X_{\bar{s}}$.
(Continuity is tested on cohomology sheaves: see section 5.)
\item[(2)]
Similarly, $P(X \times_{s} \eta, \Ql)$ denotes the category of perverse sheaves
 $\cF$ in $P(X_{\bar{s}}, \Ql)$ together with a continuous action of 
 $\Gal(\overline{\eta}/\eta)$ which is compatible with the action on $X_{\bar{s}}$.
\end{enumerate}
\end{Def}

Analogously, assuming $k(s)$ has cardinality $q$, we let $D^{b,\rm Weil}_c(X, \Ql)$ denote the category of pairs $(\cF, F^{\ast}_q)$ where $\cF$ belongs to $D^b_c(X_{\bar{s}}, \Ql)$, and $F^{\ast}_q : \Fr_q^{\ast} \cF 
~ \tilde{\rightarrow} ~ \cF$ is an isomorphism in that category.  The definition of $P_{\rm Weil}(X,\Ql)$ is similar. 

Any choice of geometric Frobenius $\Fr_q \in \Gamma$ gives rise naturally to a functor
$$ D^b_c(X \times_{s} \eta, \Ql) \lto D^{b,{\rm Weil}}_c(X, \Ql). $$
This induces an analogous functor on the categories of perverse sheaves.

We also have a natural functor
$$ P(X_s, \Ql) \lto P_{\rm Weil}(X, \Ql), \quad
    \cF_0 \mapsto (\cF, F_q^\ast), $$
which is a full embedding, and its essential image is a subcategory which
is stable under extensions and subquotients; see \cite{BBD}, Prop. 5.1.2.
In particular, a Weil perverse sheaf which has a filtration such that the graded
pieces are all defined over $s$, lies in the essential image of this
functor.

\subsection{The affine Grassmannian and flag variety}

We briefly recall the definition of the affine Grassmannian $\Grass$, and the affine flag variety $\Fl$.  


We fix a field $k$ and suppose $G$ is a split connected reductive group over $k$.  Fix a maximal torus $T$ and a Borel subgroup $B$ containing $T$.  Let $\Lambda_+ \subset X_*(T)$ denote the set of $B$-dominant integral coweights for $G$.  By $W$ we denote the finite Weyl group $N_G(T)/T$,
by $W_{\rm aff}$ the affine Weyl group, and by $\widetilde{W}$ the extended
affine Weyl group $N_G(T)/T_\mathcal O$, where $\mathcal O = k\xt$.

Consider $G(k\xT)$ as an ind-scheme
over $k$. Denote by $\cB \subset G(k\xT)$ the standard Iwahori
subgroup, i. e. the inverse image of $B$ under the projection
 $G(k\xt) \lto G(k)$.

The affine Grassmannian $\Grass$ (over the field $k$) is the quotient 
(of $fpqc$-sheaves)
 $G(k\xT)/G(k\xt)$; it is an ind-scheme. If $G=GL_n$ and $R$ is a $k$-algebra, $\Grass(R)$
is the set of all $R\xt$-lattices in $R\xT^n$. If $G=GSp_{2n}$, it is the set of
lattices in $R\xT^{2n}$ which are self-dual up to an element in $R\xt^\times$. 
In the same way, one can construct the
affine Grassmannian over $\Z$.
By the Cartan decomposition we have a stratification into $G(k\xt)$-orbits:
$$ \Grass = \coprod_{\mu \in \Lambda_+} G(k\xt)\mu G(k\xt)/G(k\xt). $$
Here we embed $X_*(T)$ into $G(k\xT)$ by the rule $\mu \mapsto \mu(t) \in T(k\xT)$.  We will denote the $G(k\xt)$-orbit of $\mu$ simply by $\cQ_\mu$ in the sequel.  The closure relations are determined by the standard partial order on dominant coweights: $\cQ_\lambda \subset \overline{\cQ}_\mu$ if and only if $\mu - \lambda$ is a sum of $B$-positive coroots.

The affine flag variety $\Fl$ over $k$ is the quotient $G(k\xT)/\cB$;
it is an ind-scheme, too. As before, we have a modular interpretation:
 e.g., for $GL_n$, $\Fl(R)$ is the space of complete lattice chains in $R\xT^n$.
The Iwahori group $\cB$ acts on $\Fl$, and we get a decomposition
$$ \Fl = \coprod_{w \in \Wext} \cB w \cB/\cB. $$
Here we embed $\widetilde{W} = X_*(T) \rtimes W$ into $G(k\xT)$ by using the aforementioned embedding of $X_*(T)$, and by choosing a representative for $w \in W = N_G(T)/T$ in the group $N_G(T) \cap G(k\xt)$.  
We write $\cB_w = \cB w \cB/\cB$; $\cB_w$ is called the Schubert cell
associated to $w$.
Its closure, which we denote by $\cB_{\overline{w}}$,
is a finite-dimensional projective variety, and is called a Schubert 
variety. We have 
$$ \cB_{\overline{w}} = \coprod_{v \le w} \cB_v. $$
Here $\le$ denotes the Bruhat order on $\widetilde{W}$ determined by the affine reflections $S_{\rm aff}$ through the walls of the {\em opposite} base alcove, i.e., the alcove $w_0(A)$, where $w_0$ is the longest element in $W$ and $A \subset X_*(T) \otimes \mathbb R$ is the alcove fixed by $\cB$ (see \cite{IM} or \cite{HKP}).

The length function $\ell: \widetilde{W} \rightarrow \mathbb Z_{\geq 0}$ used thoughout this paper is always defined with respect to the Coxeter system $(W_{\rm aff},S_{\rm aff})$.  Thus, $\cB_w$ is an affine space isomorphic to ${\mathbb A}^{\ell(w)}$.  Also, we will often write $\ell(\mu)$ instead of $\ell(t_\mu)$, where $t_\mu$ is the translation element in $\widetilde{W}$ corresponding to $\mu \in X_*(T)$.

The analogue of the unipotent radical $N^-$ of the opposite Borel $B^-$ 
in the finite-dimensional case, is the subgroup 
$\cNN \subset G(k\xT)$ which is by definition the inverse image of
 $N^-$ under the projection $G(k[t^{-1}]) \lto G(k)$. We can then decompose
 $\Fl$ into $\cNN$ orbits:
$$ \Fl = \coprod_{w \in \Wext} \cNN w\cB/\cB. $$
We write $\cB^w = \cNN w\cB/\cB$. This is an ind-scheme of 'finite 
codimension'. More precisely, $\cB^w \cap \cB_v \ne \emptyset$ if
and only if $w \le v$, and in this case the intersection has 
dimension $\ell(v) -\ell(w)$; cf. \cite{F}. We denote the closure of $\cB^w$ by
 $\cB^{\overline{w}}$; then 
$$ \cB^{\overline{w}} = \coprod_{v \ge w} \cB^v. $$

\subsection{Intersection complexes}

Let $X$ be a finite-type reduced and geometrically irreducible scheme over a finite field $k$
(or, we may assume $k$ is algebraically closed, or that $X$ is finite-type
over a discrete valuation ring with finite residue field and with fraction field $k$).  For any open immersion $j : U \rightarrow X$ of a smooth irreducible dense open subset, we define $IC(X) = j_{!*}(\Ql)$, where $j_{!*}$ is the Goresky-MacPherson-Deligne extension functor (for the middle perversity), cf. \cite{GM}, \cite{BBD}.  The shift 
$IC(X)[{\rm dim}(X)]$ is an irreducible perverse sheaf on $X$, independent of the choice of $U$.  
If $X$ is not reduced, we set $IC(X) := IC(X_{\rm red})$.

The intersection cohomology groups of $X$ are given by $IH^i(X) = H^i(X, IC(X))$ (similarly for compact supports).

We fix a choice of $\sqrt q \in \Ql$, needed to define Tate-twists on the categories $D^{b,\rm Weil}_c(X,\Ql)$ and $P_{\rm Weil}(X,\Ql)$.  It is known that $IC(X)[{\rm dim}(X)](\frac{{\rm dim}(X)}{2})$ is a self-dual pure perverse sheaf of weight zero.

For $w \in \widetilde{W}$ we have the locally closed immersion $j_w: \cB_w \hookrightarrow \Fl$.  We define the perverse sheaf $IC_w = IC(\cB_{\bar w})[\ell(w)]$.  For $\mu \in \Lambda_+$, the perverse sheaf $IC_\mu$ on $\rm Grass$ has a similar definition.  Note that
there is no Tate-twist in these definitions. This turns out to be the most
suitable normalization for our purposes.

\subsection{Deformations of the affine Grassmannian to the flag variety}

Suppose $G = GL_n$ or $GSp_{2n}$, and $S = {\rm Spec}(\Z_p)$.  There is a deformation from
 $\Grass_{\Q_p}$ to $\Fl_{\F_p}$, i.e. an ind-scheme $M$ over
 $\Z_p$ with generic fibre $\Grass_{\Q_p}$ and special fibre
 $\Fl_{\F_p}$.   The second author and Ng\^{o} \cite{HN1} defined such a 
deformation as a union
$$
M = \bigcup_{\mu \in \Lambda_+} M_\mu,
$$
where $M_\mu$ is a finite-dimensional projective $\Z_p$-scheme with generic fiber $\overline{\cQ}_{\mu,\Q_p}$ and special fiber 
$$
\coprod_{w \in {\rm Adm}(\mu)} \cB_{w,\F_p}.
$$
Here ${\rm Adm}(\mu)$ is by definition the following finite subset of $\widetilde{W}$:
$$
{\rm Adm}(\mu) = \{ x \in \widetilde{W}~ |~ x \leq 
t_\lambda, \,\, \mbox{for some $\lambda \in W\mu$} \};
$$
we refer to \cite{HN2} for further information about this subset.  
The model $M_\mu$ is defined in terms of lattice chains, so this construction works 
only for $G=GL_n$ or $GSp_{2n}$. It is conceivable, however, that other groups can be handled by using
a Pl\"ucker description of $\Grass$ and $\Fl$, cf. \cite{H-Pluecker}.

In the function field case, suppose $G$ is an arbitrary connected reductive group over $k= \F_p$ and suppose $X$ is a smooth curve over $k$.  Then a deformation ${\rm Fl}_X$ has been constructed
by Beilinson; see \cite{Gaitsgory}.  For a distinguished point $x_0 \in X$, ${\rm Fl}_X$ is an ind-scheme over $X$ whose fiber over $x \neq x_0$ is isomorphic to the product ${\rm Grass}_{k} \times G/B$, and whose fiber over $x_0$ is isomorphic to the affine flag variety $\Fl_{k}$.  This deformation also has a modular interpretation, this time in terms of $G$-bundles on $X$.  We get a deformation over $S = {\rm Spec}(k\xt)$ by base-changing with a formal neighborhood of $x_0 \in X$.

In either the function-field or the $p$-adic setting, let ${\rm Grass}$ denote the ``constant'' deformation of the affine Grassmannian, i.e., we consider the base change ${\rm Grass}_\Z \times S$.  The modular interpretations of $M$ and ${\rm Fl}_X$ allow us to define a projective morphism of ind-schemes
$$
\pi : M \rightarrow {\rm Grass}
$$
over $S$ (similarly for ${\rm Fl}_X$).  In the $GL_n$ case for example, $\pi$ sends a lattice chain ${\mathcal L}_0 \subset {\mathcal L}_1 \subset \cdots $ in $M$ to the lattice ${\mathcal L}_0$ in ${\rm Grass}$.  Often, we will use the same letter $\pi$ to denote the morphism on the special fibers, where it is just the usual projection of ind-schemes $G(k\xT)/\cB \rightarrow G(k\xT)/G(k\xt)$.

In the sequel, when we use the symbol $\Fl$ to denote the affine flag variety, we are often thinking of it as the special fiber of one of the deformations above.  Occasionally, we even write $\Fl$ for the deformation itself.  

\subsection{Nearby cycles}

Let $X$ denote a finite-type scheme over $S$.  For $\cF \in D^b_c(X_\eta,\Ql)$, we define the {\em nearby cycles sheaf} to be the object in $D^b_c(X \times_s \eta, \Ql)$ given by
$$
R\Psi^X({\cF}) = \bar{i}^*R\bar{j}_* (\cF_{\bar{\eta}}),
$$
where $\bar{i}: X_{\bar{s}} \hookrightarrow X_{\bar{S}}$ and $\bar{j} : X_{\bar{\eta}} \hookrightarrow X_{\bar{S}}$ are the closed and open immersions of the geometric special and generic fibers of $X/S$, and $\cF_{\bar{\eta}}$ is the pull-back of $\cF$ to $X_{\bar{\eta}}$.  It is known that the functor $R\Psi^X$ preserves perversity, cf. \cite{I}.

\subsection{Equivariant sheaves}

Consider the affine flag variety $\Fl$ over $\F_p$. By definition,
a (mixed) perverse sheaf on $\Fl$ is a (mixed) perverse sheaf which is supported
on some finite-dimensional part.

\begin{Def}Denote by $a \colon \cB \times \Fl \lto \Fl$ the action of the
  Iwahori group, by $i: \Fl \lto \cB\times \Fl$ the zero section,
 and by $p_2 \colon \cB \times \Fl \lto \Fl$ the projection.
A perverse sheaf $\cF$ is $\cB$-equivariant if there exists an isomorphism
 $\varphi: a^\ast \cF \lto p_2^\ast \cF$, such that

i) $\varphi$ is rigidified along the zero section: $i^\ast \varphi = \id_\cF$,

ii) $\varphi$ satisfies the cocycle condition: $(m \times \id_{\Fl})^\ast
(\varphi) = p_{23}^\ast(\varphi) \circ (\id_\cB \times a)^\ast(\varphi)$ on
$\cB \times \cB \times\Fl$, where $m \colon \cB \times \cB \lto \cB$ is the
multiplication map (cf. \cite{KW} III.15).
\end{Def}

Note that since by definition all perverse sheaves on $\Fl$ are supported on
finite-dimensional subschemes of the affine flag variety, and since the
Iwahori subgroup acts through a finite quotient on these subschemes (if
they are suitably chosen), the fact that we deal with ind-schemes does not
pose any problems at this point.

{\em Remark.} As Kiehl and Weissauer explain in \cite{KW} III.15, in our
situation the conditions i) and ii) in the definition are automatically
satisfied; more precisely, given any isomorphism $a^\ast \cF \lto p_2^\ast
\cF$, one can change it into a rigidified isomorphism, and every rigidified
isomorphism satisfies the cocycle condition. (Here we use that the 
Iwahori acts through geometrically connected quotients on suitable 
finite pieces of the affine flag variety.)

Similarly, on the affine Grassmannian we have an action of the maximal
parahoric subgroup $K=G(\mathcal O)$, where $\mathcal O=\Z_p$ in the unequal characteristic
case, and $\mathcal O=\F_p\xt$ in the function field case.

\begin{stz} 
If $\cF$ is a $K$-equivariant perverse sheaf on $\Grass$, then its sheaf of
nearby cycles $R\Psi(\cF)$ is $\cB$-equivariant.
\end{stz}

{\em Proof.} In the function field case, this is \cite{Gaitsgory} Prop. 4; 
for the unequal characteristic case see \cite{HN1} \S 4.
\qed

We denote by $P^\cB(\Fl, \Ql)$ the category of $\cB$-equivariant
perverse sheaves on $\Fl/\F_p$.  Later on (section 4) we will consider
$\cB$-equivariant perverse sheaves which have Galois structure compatible with
the group action.  We postpone the discussion of those sheaves to that section. 

We note that the intersection complexes $IC_\mu$ on ${\rm Grass}$ are $G(\mathcal O)$-equivariant, and the intersection complexes $IC_w$ on ${\Fl}$ are $\cB$-equivariant.

\subsection{The Kottwitz conjecture}

For any ${\mathcal F} \in D^b_c(\Fl \times_s \eta, \Ql)$, one can define its semi-simple trace function 
$$
x \mapsto \Tr^{ss}(\Fr_q, \cF_x),
$$
for any $x \in \Fl(\F_q)$; cf. \cite{HN1}.  Here $\Fr_q \in \Gamma$ denotes an arbitrary geometric Frobenius (the semi-simple trace is independent of the choice).  If $\cF$ is a suitably $\cB$-equivariant object (e.g., an object of $P^{\cB}(\Fl \times_s \eta,\Ql)$, cf. section 4.1), then this function defines an element in the Iwahori-Hecke algebra
$$
\mathcal H = C_c(\cB \backslash G(\F_q\xT) / \cB)
$$
of $\cB$-bi-invariant compactly-supported $\Ql$-valued functions on $G(\F_q\xT)$.  (Remark: in the sequel, we will often conflate the notions of Iwahori-Hecke algebras and affine Hecke algebras without lingering on the difference.)

In the $p$-adic setting, fix a dominant coweight $\mu$ for the group $G = GL_n$ or $GSp_{2n}$, and the associated model $M_\mu$ from \cite{HN1}.  For the nearby cycles sheaf $R\Psi := R\Psi^{M_\mu}(IC_\mu)$, the Kottwitz conjecture is the following equality of functions in the Iwahori-Hecke algebra for $G$:
$$
\Tr^{ss}(\Fr_q, R\Psi) = \varepsilon_\mu q_\mu^{1/2} \sum_{\lambda \leq \mu} m_\mu(\lambda) z_\lambda.
$$
Here, $\varepsilon_\mu = (-1)^{\ell(\mu)}$, $q_\mu = q^{\ell(\mu)}$, $\lambda$ ranges over dominant 
coweights preceding $\mu$ in the usual partial order, $m_\mu(\lambda)$ is the multiplicity of $\lambda$ in the character of the dual group having highest weight $\mu$, and $z_\lambda$ is the Bernstein function attached to $\lambda$.  More precisely $z_\lambda = \sum_{\lambda \in W\mu} \Theta_\lambda$, where $\Theta_\lambda$ is the function $\tilde{T}_{t_{\lambda_1}} \tilde{T}^{-1}_{t_{\lambda_2}}$;  here $\lambda = \lambda_1 - \lambda_2 $ and $\lambda_i$ is dominant, and $\tilde{T}_x := q_x^{-1/2}T_x$ for $x \in \widetilde{W}$.  The symbol $T_x$ denotes the generator ${\rm char}(\cB x \cB)$ in the Iwahori-Hecke algebra.

This formula was proved in \cite{HN1} in the situation at hand.  In the function field case, the analogous formula holds.  In that analogue, $IC_\mu$ is replaced with the exterior product $IC_\mu \boxtimes \delta$ on ${\rm Grass} \times G/B$, where $\delta$ is the skyscraper sheaf supported at the base point of $G/B$.  This is a consequence of Gaitsgory's work \cite{Gaitsgory}, and holds for any connected reductive group $G$.

We make extensive use of this formula, in particular in the computation of the examples in section 9, and in the proof that $R\Psi$ satisfies the property (P) (cf. section 6).

\section{Inverse Kazhdan-Lusztig polynomials}

\subsection{The sheaf $IC(\cB^{\overline{v}})$}

Let us recall some of the notation introduced above.  Let $k = \bar{\F}_p$.
We denote by $\cB$ the Iwahori subgroup of $G(k\xT)$, i.e. the inverse
image of the fixed Borel group $B$ under the projection $G(k\xt) \lto G$.
By $\cNN$ we denote the inverse image of the unipotent radical of the
opposite Borel under the projection $G(k[t^{-1}]) \lto G$. Furthermore,
let $\cNN(n)$ be the inverse image of $T(k[t^{-1}]/t^{-n})$ under
the projection $\cNN \lto G(k[t^{-1}]/t^{-n})$ (cf. \cite{F}).  Note that ${\rm Lie}(\cN^-)$ is generated by the negative affine roots (recall we embed $X_*(T) \hookrightarrow T(k((t)))$ by $\lambda \mapsto \lambda(t)$).

For elements $w$ and $y$ in the extended affine Weyl group $\widetilde{W}$,
we denote by $\cB_y$ the $\cB$-orbit of $y$, by $\cB^w$ the $\cNN$-orbit
of $w$, by $\cB^{\overline{w}}$ its closure, by $\cB_y^w$ the intersection
$\cB_y \cap \cB^w$, and by $\cB_y^{\overline{w}}$ the intersection
$\cB_y \cap \cB^{\overline{w}}$.  Recall the basic fact (cf. \cite{F}) that $\cB^w_y \neq \emptyset \Leftrightarrow w \leq y$ in the Bruhat order on $\widetilde{W}$.

Further, we define the opposite Iwahori subgroup $\cB^-$ to be the inverse image of the opposite Borel subgroup $B^-$ under the projection $G(k[t^{-1}]) \rightarrow G(k)$.  Note that $\cN^- \vartriangleleft \cB^-$ and 
$\cB^- = \cN^- \cdot (\cB^- \cap T(k((t))))$.  Hence for any $w \in \widetilde{W}$, $\cN^-w* = \cB^-w*$, where $*$ denotes the base point in $\Fl$.

Fix $w \in \Wext$ such that $\Adm(\mu) \subseteq \{ v ~|~ v \le w \}$.

Let $\Omega = \bigcup_{ v \le w } v \cNN \ast$. 
This is an open subset of $\Fl$.

\begin{lem}
We have $\Omega = \coprod_{v \le w} \cB^v$.
Thus for $y \le w$ we have $\cB_y \subseteq \Omega$.
\end{lem}

{\em Proof.}
We define 
\begin{eqnarray*}
\cN^-_w & = & \cN^- \cap w \cN^- w^{-1} \\
\cN^+_w & = & \cN^+ \cap w \cN^- w^{-1}.
\end{eqnarray*}
Here $\cN^+ $ denotes the preimage of $N$ under the projection $G(k[[t]]) \rightarrow G(k)$; ${\rm Lie}(\cN^+)$ is generated by the positive affine roots.

We then have 
$$ \cN^- \cong \cN^-_w \times \prod_{\gfrac{a < 0}{w^{-1}(\alpha) > 0}}
U_\alpha,$$
and
$$ \cN^+_w = \prod_{\gfrac{\alpha>0}{w^{-1}(\alpha)<0}} U_\alpha
 = w\left( \prod_{\gfrac{\alpha<0}{w(\alpha) > 0}} U_\alpha \right)
 w^{-1}. $$
Furthermore, it is easy to see that
$\cB_w = \cN^+_w w \ast$ and $\cB^w = \cN^-_w w \ast$.

From the above we have
$$ v \cN^-  = v \cN^-_{v^{-1}} \cdot (v^{-1} \cN^+_v v) = \cN^-_v \cN^+_v v. $$
Thus clearly
$$ \cB^v = \cN^-_v v \ast \subseteq v \cN^- \ast, $$
so
$$ \coprod_{v \le w} \cB^v \subseteq \bigcup_{v \le w} v \cN^- \ast.$$
On the other hand, $\cN^+_v v = \cB_v \subseteq \coprod_{v' \le v}
\cB^{v'}$, hence $\cN^-_v \cN^+_v v \subseteq \coprod_{v' \le v} \cB^{v'}$
and thus 
$$ \bigcup_{v \le w} v \cN^- \ast \subseteq \coprod_{v \le w} \cB^v.$$
\qed

For $n$ sufficiently large, the quotient 
$\cNN(n) \backslash \Omega$ exists; cf.  \cite{F} (choose $n$ large enough so that $\cN^-(n) \subset 
\cN^-_v$, for every $v \le w$). The reason 
is that we have a product decomposition
$$ v\cNN \ast = \cNN_v \cN^+_v v \ast \cong \cNN(n) \times \A^a \times \A^b, $$
where $b = {\rm dim}(\cN^+_v) = \ell(v)$, and $a$ is the number of affine roots in $\cN^-_v$ which are not 
in $\cN^-(n)$; cf. \cite{F}. Thus we see that more precisely, the quotient
 $\cNN(n) \backslash \Omega$ can be covered by affine spaces,
and in particular is smooth.

In the following discussion, we often tacitly identify $\cB_y = \cN^+_y y$, $\cB^y = 
\cN^-_y y$ and $\cB^y \times \cB_y = \cN^-_y\cN^+_y y = y\cN^- \ast$.  We have a commutative diagram $(v \le y \le w$) 
$$ \xymatrix{
\cB_y^{\overline{v}} \ar[d]^{\cong} \ar@{^{(}->}[r] & \cB_y \ar[d]^{\cong}
\ar@{^{(}->}[r] & \cB^y \times \cB_y \ar[d] \ar@{^{(}->}[r]  & \Omega
\ar[d]^{\pi} \\
\pi(\cB_y^{\overline{v}}) \ar@{^{(}->}[r] & \pi(\cB_y) \ar@{^{(}->}[r] & 
        \cNN(n) \backslash \cB^y \times \cB_y \ar@{^{(}->}[r] & \cNN(n)
        \backslash \Omega.
} $$

Note that
 $$ \pi(\cB_y^{\overline{v}} \cap \Omega) 
  = \pi(\cB_y) \cap \pi(\cB^{\overline{v}} \cap \Omega). $$

Furthermore $\pi(\cB^v \cap  \Omega)$ is isomorphic to an affine space
(thus in particular is smooth),
and its codimension in $\cNN(n) \backslash \Omega$ is $\ell(v)$.
Finally, we have
 $$ \cNN(n) \backslash (\cB^y \times \cB_y) = 
    (\cNN(n) \backslash \cB^y) \times \cB_y =: \cA_y $$
and this is an open subset of $\cNN(n) \backslash \Omega$,
and $\cA_y \cap \pi(\cB^{\overline{v}} \cap \Omega) = 
     (\cNN(n)\backslash \cB^y) \times \cB_y^{\overline{v}}$ 
is an open subset of
 $\pi(\cB^{\overline{v}} \cap \Omega)$.

By taking the quotient of the morphism $\cNN \times \Omega \longrightarrow
\Omega \times \Omega$, $(g, x) \mapsto (x, gx)$ giving the action of $\cNN$ on
$\Omega$, we get the groupoid
$$ \cNN(n) \backslash \cNN \times^{\cNN(n)} \Omega \longrightarrow \cNN(n)
\backslash \Omega \times \cNN(n) \backslash \Omega $$
(cf. \cite{F}).

\begin{stz}
This is a smooth groupoid, i.e. the maps $\cNN(n) \backslash \cNN \times^{\cNN(n)} \Omega \longrightarrow \cNN(n)
\backslash \Omega$, $(g,x) \mapsto x$, and $\cNN(n) \backslash \cNN \times^{\cNN(n)} \Omega \longrightarrow \cNN(n)
\backslash \Omega$, $(g,x) \mapsto gx$, are smooth.
\end{stz}


{\em Proof.}
Consider the maps 
$$ \phi_1, \phi_2: \cNN \times \Omega \lto \Omega, \quad 
   (g,x) \mapsto x, \text{ resp. } (g,x) \mapsto gx. $$
These are formally smooth maps since $\cNN$ is 
formally smooth. 
They induce maps
$$ \psi_1, \psi_2: \cNN(n) \backslash \cNN \times^{\cNN(n)} \Omega
   \lto \cNN(n) \backslash \Omega, $$
where $\psi_1(g, x) = x$, $\psi_2(g,x) = gx$.

We must show that $\psi_1$ and $\psi_2$ are smooth.
We have a commutative diagram
$$ \xymatrix{
\cNN \times \Omega \ar[r]^{\phi,\text{fs}} \ar[d]^{\text{fs}} & \Omega \ar[d]^{\text{fs}} \\
\cNN(n)\backslash \cNN \times^{\cNN(n)} \Omega \ar[r]^<<<<<\psi & 
\cNN(n) \backslash \Omega
} $$
where $\phi$ is either $\phi_1$ or $\phi_2$, $\psi $ is the corresponding
map and 'fs' means that
we already know that this morphism is formally smooth.

Since $\cNN(n)\backslash \cNN \times^{\cNN(n)} \Omega \lto 
\cNN(n) \backslash \Omega$ is a morphism of finite type
between Noetherian schemes (over $\overline{\F}_p$), in order
to prove that it is smooth, it is sufficient to show the 
infinitesimal lifting criterion for local Artin rings with residue
class field $\overline{\F}_p$. (see \cite{SGA1}, Exp. III Thm. 3.1).

So let $A$ be a local Artin ring with residue field $\overline{\F}_p$,
let $I \subset A$ be a nilpotent ideal, and look at a diagram
$$ \xymatrix{
\Spec A/I \ar[r] \ar[d] & \cNN(n)\backslash \cNN \times^{\cNN(n)} \Omega \ar[d] \\
\Spec A \ar[r] & \cNN(n) \backslash \Omega.
} $$ 
We have to find a morphism $\Spec A \lto \cNN(n)\backslash 
\cNN \times^{\cNN(n)} \Omega$ which makes the diagram commutative.

Since the morphism $\cNN \times \Omega \lto \cNN(n)\backslash 
\cNN \times^{\cNN(n)} \Omega$ is surjective on $\overline{\F}_p$-valued
points, we can extend the previous diagram to a diagram
$$ \xymatrix{
\Spec \overline{\F}_p \ar[r] \ar[d] & \cNN \times \Omega \ar[d]^{\xi}\\
\Spec A/I \ar[r] \ar[d] & \cNN(n)\backslash \cNN \times^{\cNN(n)} \Omega  \ar[d]^{\psi} \\
\Spec A \ar[r] & \cNN(n) \backslash \Omega.
} $$ 

Now since $\xi$ is formally smooth, we can insert a map $\Spec A/I 
\lto  \cNN \times \Omega$ into the diagram. Then we get a map
$\Spec A \lto \cNN \times \Omega$, since the composition $\psi\xi$
is formally smooth. Now we can define the map we are looking for as 
the composition of this map with $\xi$.

Thus we see that $\psi_1$ and $\psi_2$ are indeed smooth.
\qed

From the Proposition, we immediately get the following corollary which will be
used in the proof of the theorem of Kazhdan and Lusztig.

\begin{kor} \label{equisingularity} 
Let $v \le z \le w$.
For all $z' \in \cB^z(\F_q)$,
 $$ IC_{z'} (\pi(\cB^{\overline{v}} \cap \Omega)) \cong
   IC_{z} (\pi(\cB^{\overline{v}} \cap \Omega)). $$
\end{kor}

{\em Proof.}
By pull-back the maps $\psi_1$ and $\psi_2$ in the Proposition induce 
smooth maps
$$ \psi_1, \psi_2 : \cNN(n) \backslash \cNN \times^{\cNN(n)} 
   (\cB^{\overline{v}} \cap \Omega) \lto \pi(\cB^{\overline{v}}\cap \Omega), $$
where $\psi_1(g, x) = x$, $\psi_2(g,x) = gx$.

Now choose $g \in \cNN$ such that $gz = z'$. 
The point $(g,z)$ maps to $z$ resp. $z'$
under $\psi_1$ resp. $\psi_2$. Since both these maps are smooth,
the corollary follows. 
\qed

However, we can even say more: The sheaf $IC(\pi(\cB^{\overline{v}} \cap
\Omega))$ on $\pi(\cB^{\overline{v}} \cap \Omega) \subseteq \cNN(n) \backslash
\Omega$ is equivariant under the above groupoid (our definition of
equivariance obviously makes sense in this context). Namely, we have
$$ \psi_1^* IC(\pi(\cB^{\overline{v}} \cap \Omega)) \cong 
   IC(\cNN(n)\backslash \cNN \times^{\cNN(n)} (\cB^{\overline{v}} \cap \Omega)) \cong
   \psi_2^* IC(\pi(\cB^{\overline{v}} \cap \Omega)),$$
since $\psi_1$ and $\psi_2$ are both smooth. We can choose an isomorphism
which is rigidified along the zero section, and since $\cNN(n)\backslash \cNN$
is connected such an isomorphism will automatically satisfy the cocycle
condition.

Thus the pull-back of the sheaf $IC(\pi(\cB^{\overline{v}} \cap \Omega))$
to $\Omega$ is equivariant with respect to the action of $\cNN$.

We can now define the intersection complex $IC(\cB^{\overline{v}})$
of the anti-cell $\cB^{\overline{v}}$. Since $\cB^{\overline{v}}$ is just an
ind-scheme, it is not a priori clear how to make sense of that. However, by
looking at the quotients $\cNN(n) \backslash \Omega$, we can define a sheaf
on $\Fl$ which merits the name $IC(\cB^{\overline{v}})$.

The pull-back of the sheaf $IC(\pi(\cB^{\overline{v}} \cap \Omega))$
to $\Omega$ is independent of $n$: in fact, let $n' > n$,
and consider the diagram
$$ \xymatrix{
\Omega \ar[rr]^\pi \ar[dr]^{\pi'} & & \cNN(n) \backslash\Omega \ar[dl]_{\psi} \\
 & \cNN(n')\backslash\Omega. & \\
} $$
The map $\psi$  is an $\A^N$-bundle (for suitable $N$), and $\pi(\cB^{\overline{v}} \cap \Omega)
 = \psi^{-1}(\pi'(\cB^{\overline{v}} \cap \Omega))$. Thus the pull-back
of $IC(\pi'(\cB^{\overline{v}} \cap \Omega))$ along $\psi$ is just
 $IC(\pi(\cB^{\overline{v}} \cap \Omega))$. This means that the pull-back
 of $IC(\pi(\cB^{\overline{v}} \cap \Omega))$ to $\Omega$ is independent of
 $n$. Since it is obviously compatible with enlarging $\Omega$, we indeed
 get a sheaf on $\Fl$, which we denote by $IC(\cB^{\overline{v}})$.

We will only ever need to work with the restriction of
$IC(\cB^{\overline{v}})$ to finite-dimensional subschemes of $\Fl$. In
particular we do not really need to define this sheaf on $\Fl$ --- we could
always work on a suitable quotient $\cNN(n) \backslash \Omega$. 

We will need the following property of $IC(\cB^{\overline{v}})$. When
restricted to an intersection $\cB^{\overline{v}}_y$, we get the
  intersection complex of that space:
$$ IC(\cB^{\overline{v}})|_{\cB^{\overline{v}}_y} = IC
(\cB^{\overline{v}}_y).$$
This follows from 

\begin{stz} \label{restriction} We have
$$ IC(\pi(\cB^{\overline{v}} \cap \Omega))|_{\pi(\cB^{\overline{v}}_y)}
   = IC (\cB^{\overline{v}}_y). $$
\end{stz}

{\em Proof.} This follows immediately from the fact that 
$\cA_y \cap \pi(\cB^{\overline{v}} \cap \Omega) = 
     (\cNN(n)\backslash \cB^y) \times \cB_y^{\overline{v}}$ 
is an open subset of
 $\pi(\cB^{\overline{v}} \cap \Omega)$ and that 
 $\cNN(n)\backslash \cB^y$ is just an affine space. \qed

\subsection{The theorem of Kazhdan and Lusztig}

We start with two lemmas, and then prove the theorem of Kazhdan and Lusztig.

\begin{lem} \label{numberofpoints} 
We have
$$ \# \cB^z_y (\F_q) = R_{z,y}(q). $$
\end{lem}

{\em Proof.} See \cite{H}, section 2. \qed

\begin{lem} \label{gmaction} On $\cA_y$, which is isomorphic to an
affine space ${\mathbb A}^N$, we have a $\G_m$-action given by
$$ \lambda \cdot (z_1, \dots, z_N) = (\lambda^{a_1} z_1, \dots, \lambda^{a_N} z_N), $$
for certain $a_i > 0$, which preserves 
 $\pi(\cB^{\overline{v}} \cap \Omega) \cap \cA_y$ and
$\pi(\cB^{\overline{v}} \cap \cB_y)$.
\end{lem}

{\em Proof.} We decompose $\cA_y = \cNN(n)\backslash \cB^y \times \cB_y$,
and we will define $\G_m$-actions on the two factors separately.
With respect to this decomposition, we have 
 $\pi(\cB^{\overline{v}} \cap \Omega) \cap \cA_y = \cNN(n)\backslash \cB^y
 \times \cB^{\overline{v}} \cap \cB_y$
and $\pi(\cB^{\overline{v}} \cap \cB_y) = \{y\} \times \cB^{\overline{v}}
\cap \cB_y$.
Thus the $\G_m$-action on the first factor must just be a contracting
$\G_m$-action which fixes the origin.

To define the action on the second factor, consider the action
of the torus $T(k)$ on $\cB_y$. Obviously it fixes
$\cB^{\overline{v}} \cap \cB_y$, so we only have to find
a cocharacter $\phi: \G_m \lto T$ such that the induced $\G_m$-action
is contracting.

This is possible because the set of finite roots $\alpha$
such that for some $n$, the affine root $(\alpha, n)$ 'occurs' in $\cB_y$
is the union of the sets
$$
\{ \alpha > 0 ~ | ~ \langle \alpha,\lambda \rangle \le 0, \,\, \mbox{and $\tilde{y}\alpha < 0$ if $\langle \alpha,\lambda \rangle = 0$} \}
$$
and 
$$\{ \alpha < 0 ~ | ~ \langle \alpha, \lambda \rangle \le -1, \,\, \mbox{and $\tilde{y}\alpha < 0$ if $\langle \alpha, \lambda \rangle = -1$} \},
$$
where $y^{-1} = \tilde{y}t_\lambda$, $\tilde{y} \in W$, $\lambda \in X_\ast(T)$.  We may put $\phi = 
-m\lambda - \tilde{y}^{-1}\rho^\vee$ for any sufficiently large integer $m$.
\qed

\begin{thm}
\label{invKLpol}
 {\rm (Interpretation of inverse Kazhdan-Lusztig polynomials. 
(\cite{KL2}, Prop. 5.7))}

Let $v \le y$. Then $IH^i(\cB_y^{\overline{v}})=0$ for $i$ odd,
and 
$$ Q_{v,y}(q)  =  
     \sum_{i\ge 0} \dim IH^{2i}(\cB_y^{\overline{v}})q^i 
 = \Tr(\Frob_q, IH^\bullet(\cB_y^{\overline{v}})).
$$
\end{thm}

{\em Remark.} This result is stated in \cite{KL2} as a proposition
which 'should not be difficult to prove'.

Kashiwara and Tanisaki give a proof of the analog of this result
in the setting of mixed Hodge modules in \cite{KT}, Thm. 6.6.4; 
in particular they work over the complex numbers. 

{\em Proof.} Given the lemmas above, the theorem is proved much in
the same way as Theorem 4.2 in \cite{KL2} which gives an analogous
interpretation for the Kazhdan-Lusztig polynomials in the finite 
dimensional case.

Choose $\Omega = \bigcup_{ v \le w } v \cNN \ast$ 
sufficiently large, and denote the projection 
 $\Omega \lto \cNN(n) \backslash \Omega$ by $\pi$ as before. 

We write $\xB^?_?$ as an abbreviation for $\pi(\cB^?_? \cap \Omega)$.
Note that $\xB^v_y \cong \cB^v_y$.
Denote by $\sH^i(\xB^{\overline{v}})$ the $i$-th cohomology sheaf 
of $IC(\xB^{\overline{v}})$, and by $\sH^i_z(\xB^{\overline{v}})$ its
stalk in a point $z$.

First of all, Lemma 4.5a) in \cite{KL2} together with Lemma 
\ref{gmaction}, the fact that $\xB^y \times \xB^{\overline{v}}_y$
is open in $\xB^{\overline{v}}$, and the fact that $\xB^y$ is just an affine space,
yield that 
 $IH^i(\xB^{\overline{v}}_y) \cong \sH^i_y(\xB^{\overline{v}})$.
In the sequel we will thus always work with 
 $\sH^i_y(\xB^{\overline{v}})$.

For $z \ge v$ consider the property

{\bf P($z$):} For any $z' \in \xB^z(\F_{q})$, we have
 $\sH^i_{z'}(\xB^{\overline{v}})=0$ for odd $i$, and for even $i$
 all eigenvalues of $\Frob_q$ on $\sH^i_{z'}(\xB^{\overline{v}})$
 are equal to $q^{i/2}$.

Obviously {\bf P($v$)} is true. Now fix $y > v$ and assume that 
{\bf P($z$)} holds
for all $v \le z < y$. We want to show that this implies 
{\bf P($y$)}.

By assumption $(\xB^{\overline{v}} \cap \cA_y) - \xB^y$ is 
very pure. Since
$$ (\xB^{\overline{v}} \cap \cA_y) - \xB^y \cong \xB^y \times 
 (\xB^{\overline{v}}_y - \{ y \}) $$
and $\xB^y$ is smooth, we see that 
 $\xB^{\overline{v}}_y - \{ y \}$
is very pure. Thus by Lemma \ref{gmaction} and Lemma 4.5b) in \cite{KL2},
 $\xB^{\overline{v}}_y$ is very pure. This implies that
 $\xB^{\overline{v}} \cap \cA_y$ is very pure.

Now the Lefschetz trace formula implies (use Lemma \ref{numberofpoints} 
and Lemma \ref{equisingularity}) 
\begin{eqnarray*}
&& \Tr(\Fr_q, IH^\bullet_c(\xB^{\overline{v}} \cap \cA_y)) \\
& = & \sum_{v \le z \le y} \sum_{z' \in (\xB^z \cap \cA_y)(\F_q)}
        \Tr(\Fr_q, \sH^\bullet_{z'}(\xB^{\overline{v}})) \\
& = & \sum_{v \le z \le y } q^{\dim \xB^y} R_{z,y}(q) 
        \Tr(\Fr_q, \sH^\bullet_z(\xB^{\overline{v}})).
\end{eqnarray*}

By Poincar\'e duality we have
$$ \Tr(\Fr_q, IH^\bullet_c(\xB^{\overline{v}} \cap \cA_y)) = 
   q^{\dim (\xB^{\overline{v}} \cap \cA_y)} 
    \Tr(\Fr_q^{-1}, IH^\bullet(\xB^{\overline{v}} \cap \cA_y)) , $$
and \cite{KL2}, Lemma 4.5 gives us 
$$ \Tr(\Fr^{-1}_q, IH^\bullet(\xB^{\overline{v}} \cap \cA_y)) =
   \Tr(\Fr^{-1}_q, \sH^\bullet_y(\xB^{\overline{v}})).   $$

Altogether, we now get
\begin{equation} \label{recursive_characterization}
 q^{\ell(y) - \ell(v)} \Tr(\Fr_q^{-1}, \sH^\bullet_y(\xB^{\overline{v}})) = 
   \sum_{ v \le z \le y} R_{z,y} \Tr(\Fr_q,
   \sH^\bullet_z(\xB^{\overline{v}})).
\end{equation}

This equation is completely analogous to (4.6.4) in \cite{KL2},
and exactly the same arguments as in loc. cit. now yield that, $\forall i$,
$$ \Tr(\Fr_q, \sH^i_y(\xB^{\overline{v}})) \in \Z[q], $$
which implies that {\bf P($y$)} holds. This proves the second
equality in the statement of the theorem.

Furthermore, since the inverse Kazhdan-Lusztig polynomials
are characterized by
\begin{gather} Q_{v,v} =1, \quad q^{\ell(y) - \ell(v)} Q_{v,y}(q^{-1}) =
                     \sum_{v \le z \le y} R_{z,y} Q_{v,z}, \notag\\
 \deg Q_{v,y} \le (l(y)-l(v)-1)/2, \,\,\, (v < y)\notag
\end{gather}
we also get the first equality in the statement of the theorem 
from equation (\ref{recursive_characterization}). \qed

\section{The main theorem}

\subsection{Decomposition of $\cB$-equivariant sheaves}

We begin with a discussion of the categories of perverse sheaves to which our results apply.  We will denote by 
$P(\Fl \times_s \eta, \Ql)$ the full subcategory of $D^b_c(\Fl \times_s \eta, \Ql)$ consisting of middle perverse sheaves endowed with a compatible action of the Galois group $\Gamma = {\rm Gal}(\bar{\eta}/\eta)$.  (For the precise meaning of ``compatible'', see section 5.1).  We let $P^{\mathcal B}(\Fl \times_s \eta, \Ql)$ denote the full subcategory of ${\mathcal B}$-equivariant objects of $P(\Fl \times_s \eta, \Ql)$: these are perverse sheaves on $\Fl_{\bar s}$ which are equivariant for $\Gamma$ and ${\mathcal B}$, in a compatible way.

Similarly, we let $P_{\rm Weil}(\Fl,\Ql)$ denote the category of pairs $({\mathcal F}, F^*_q)$ where ${\mathcal F}$ is a perverse sheaf on $\Fl_{\bar{s}}$ and $F^*_q : \Fr^*_q{\mathcal F} \tilde{\rightarrow} {\mathcal F}$ is an isomorphism of perverse sheaves (a Weil structure on ${\mathcal F}$).  Further, let $P^{\mathcal B}_{\rm Weil}(\Fl,\Ql)$ denote the full subcategory of ${\mathcal B}$-equivariant objects of $P_{\rm Weil}(\Fl,\Ql)$:  the 
${\mathcal B}$-equivariance of ${\mathcal F}$ is assumed compatible with $F^*_q$.

\begin{lem}
The subcategory $P^{\mathcal B}_{\rm Weil}(\Fl,\Ql)$ of $P(\Fl,\Ql)$ is stable under formation of kernels and cokernels.  Objects in this subcategory have finite length.  The irreducible objects are intersection complexes of the form $IC({\mathcal L})[\ell(w)] = j_{!*}{\mathcal L}[\ell(w)]$, where $j: {\mathcal B}_w \rightarrow \Fl$ is the locally closed immersion of an Iwahori-orbit ${\mathcal B}_w$ in $\Fl$, and ${\mathcal L}$ is a ${\mathcal B}$-equivariant (thus geometrically constant) irreducible Weil sheaf on ${\mathcal B}_w$. 
\end{lem}

{\em Proof.}  The statements concerning kernels and cokernels, and finite length, are clear.  We prove that all the irreducible subquotients of ${\mathcal F} \in P^{\mathcal B}_{\rm Weil}(\Fl,\Ql)$ taken in this category are of the required form $IC({\mathcal L})[\ell(w)]$, by Noetherian induction on ${\rm supp}({\mathcal F})$.  Indeed, we may find an open immersion $j: {\mathcal B}_w \rightarrow {\rm supp}({\mathcal F})$ such that $j^*{\mathcal F}$ is a (necessarily ${\mathcal B}$-equivariant) Weil-sheaf on ${\mathcal B}_w$ (up to a shift by $\ell(w)$).  Now the kernel and cokernel of the adjunction map $j_!j^*{\mathcal F} \rightarrow {\mathcal F}$ belong to $P^{\mathcal B}_{\rm Weil}(\Fl,\Ql)$ and are supported on proper closed ${\mathcal B}$-invariant subsets of ${\rm supp}({\mathcal F})$.  Further, the kernel of the canonical surjection $j_!j^*{\mathcal F} \rightarrow j_{!*}j^*{\mathcal F}$ also has this property.  Since for the action of ${\mathcal B}$ on the orbit $\cB_w$, the stabilizer of 
any point is geometrically connected, any ${\mathcal B}$-equivariant Weil sheaf on ${\mathcal B}_w$ is geometrically just a constant sheaf.

The fact that the 
subquotients of ${\mathcal F}$ have the required form is a consequence of these remarks.  (Note that $j_!$ etc. denotes the {\em derived} functor here.  We have used the fact that the morphism $j$ is affine to ensure that $j_!$ preserves perversity and that $j_{!*}$ is a quotient of $j_!$.)  

It also follows from this that the irreducible objects of $P^{\mathcal B}_{\rm Weil}(\Fl,\Ql)$ are of the stated form. 
\qed

\begin{kor}
Objects in $P^{\cB}_{\rm Weil}(\Fl,\Ql)$ are mixed.
\end{kor}

Now consider a Weil-perverse sheaf $\cF \in P^\cB_{\rm Weil}(\Fl, \Ql)$ such that
for all $x \in \Fl(\F_q)$, we have $\Tr(\Fr_q, \cF_x) \in \Z[q^{1/2}, q^{-1/2}]$.
We will denote the full subcategory of those sheaves by $P^\cB_{\sqrt q}(\Fl, \Ql)$.  Similarly, let 
$P^\cB_q(\Fl,\Ql)$ denote the full subcategory of sheaves such that 
$\Tr(\Fr_q, \cF_x) \in \Z[q,q^{-1}]$ for all $x \in \Fl(\F_q)$.  The category $P^\cB_{\sqrt q}(\Fl, \Ql)$ gives rise via the sheaf-function dictionary precisely to elements in the Iwahori-Hecke algebra (for this reason it is sometimes called the Hecke category).

\begin{lem}
The subcategory $P^\cB_{q}(\Fl, \Ql)$ in $P^\cB_{\rm Weil}(\Fl, \Ql)$ is 
stable under extensions, quotients, and subquotients, and under Verdier duality.  The irreducible objects of $P^\cB_q(\Fl,\Ql)$ are of the form $IC_w(i)$, for $w \in \widetilde{W}$, $i \in \Z$.
\end{lem}

{\em Proof.} 
The statements concerning extensions, quotients, and Verdier duality are clear.  Let ${\mathcal F} \in P^\cB_q(\Fl,\Ql)$.  We prove all of its irreducible subquotients in the category $P^\cB_{\rm Weil}(\Fl,\Ql)$ are of the form $IC_w(i)$, by Noetherian induction on ${\rm supp}({\mathcal F})$. 
Replacing ${\mathcal F}$ by its semi-simplification, we may assume ${\mathcal F}$ is semi-simple.
Choose an open immersion $j: {\mathcal B}_w \rightarrow {\rm supp}({\mathcal F})$, such that $j^*{\mathcal F}[-\ell(w)]$ is a (geometrically constant) Weil sheaf ${\mathcal L}$.  By hypothesis, $\Tr(\Fr_q,{\mathcal L}_w) \in \Z[q,q^{-1}]$.  
The linear independence of characters shows that the eigenvalues of Frobenius on ${\mathcal L}_w$ are of form $q^i$, for various integers $i$, in other words, ${\mathcal L}$ is a sum of Weil sheaves on ${\mathcal B}_w$ of the form $\Ql(i)$, $i \in \Z$.  By the previous Lemma, the intersection complex $IC({\mathcal L})[\ell(w)]$ is a subquotient of ${\mathcal F}$ in the category $P^\cB_{\rm Weil}(\Fl,\Ql)$.  The quotient ${\mathcal F}/IC({\mathcal L}[\ell(w)])$ belongs to 
$P^\cB_q(\Fl,\Ql)$ (because by Kazhdan-Lusztig \cite{KL2}, $IC({\mathcal L}[\ell(w)]) \in P^\cB_q(\Fl,\Ql)$), and has support a proper closed subset of ${\rm supp}({\mathcal F})$.  Applying the induction hypothesis to this quotient proves that all the irreducible subquotients of ${\mathcal F}$ have the required form.  

The fact that the subcategory $P^\cB_q(\Fl,\Ql)$ of $P^\cB_{\rm Weil}(\Fl,\Ql)$ is closed under formation of subobjects follows.
\qed
 
\begin{kor} 
Let $\cF$ be an object of $P^\cB_{q}(\Fl, \Ql)$. 
Then the constituents of the Jordan-H\"older series of $\cF$ have 
the form $IC_w(-i)$, $w \in \Wext$, $i\in \Z$.
Each of them appears with a certain multiplicity, which we denote by
 $m(\cF, w,i)$.
\end{kor}

\subsection{Relation to the Hecke algebra}
\label{hecke_algebra}

Denote by $\cH = \bigoplus_{w\in \Wext} \Z[q^{-1/2}, q^{1/2}] T_w$ the 
Iwahori-Hecke algebra associated to $G$.

As usual, we write $q_w = q^{\ell(w)}$, $\varepsilon_w = (-1)^{\ell(w)}$.

A perverse sheaf $\cF \in P^\cB_{\sqrt q}(\Fl, \Ql)$ gives rise to a function 
 $\Tr(\Fr_q, \cF) = \sum_w \Tr(\Fr_q, \cF_w) T_w \in \cH$, and we have a
 sheaf-function dictionary \`a la Grothendieck.

On the other hand, Kazhdan and Lusztig have defined a different basis $C_w$
of $\cH$. We slightly change the normalization, and work with elements
$C''_w$ instead. They are related to the $C_w$ resp. to the elements $C'_w$
sometimes appearing in the literature as follows. By definition, $C''_w$ is
the function $\Tr(\Fr_q, IC({\mathcal B}_w)[\ell(w)])$. Also by definition,  
 $C'_w = \Tr(\Fr_q, IC({\mathcal B}_w)(\ell(w)/2))$, so we have $C''_w = \varepsilon_w
q_w^{1/2} C'_w$. But we also have $\sigma(C_w) = \varepsilon_w C'_w$,
where  $\sigma$ is the involution of the Hecke algebra given by 
 $q \mapsto q^{-1}$, $T_w \mapsto \varepsilon_w q_w^{-1} T_w$. Thus
 $C''_w = q_w^{1/2} \sigma(C_w)$.
The base change matrix giving the relation between the $T_w$ and the
$C''_w$ is given by Kazhdan-Lusztig polynomials resp. inverse
Kazhdan-Lusztig polynomials:
\begin{eqnarray*} T_w & = & \sum_x \varepsilon_w Q_{x,w} C''_x, \\
C''_w & = & \sum_x \varepsilon_w P_{x,w} T_x.
\end{eqnarray*}

Now if $\cF \in P^\cB_q(\Fl, \Ql)$, the multiplicities $m(\cF, w, i)$
defined above are just the coefficients of $\Tr(\Fr_q, \cF)$ with respect
to the basis $C''_w$:
$$  \Tr(\Fr_q, \cF) = \sum_w \left( \sum_i m(\cF,w,i)q^i \right) C''_w. $$
Thus, knowing the function for $\cF$ explicitly is equivalent to knowing all the multiplicities $m(\cF,w,i)$ for $\cF$.

For later use note that the identity $\overline{C''_w} = q_w^{-1}C''_w$, together with the fact that Verdier duality descends to the Kazhdan-Lusztig involution on the Iwahori-Hecke algebra, gives the following identity.

\begin{lem} \label{proto_palindromic}
We write $m(\cF,w) = \sum_i m(\cF,w,i)q^i$ for $\cF \in P^{\cB}_{q}(\Fl,\Ql)$.  Then
$$
m(D\cF,w) = q^{-1}_w \overline{m(\cF,w)}.
$$
\end{lem}

\subsection{The property (P)}

We say that a sheaf $\cF \in P^\cB_{\sqrt q}(\Fl, \Ql)$ satisfies the 
property (P) (for ``palindromic''), if for all $y \in \Fl(\F_q)$ we have
$$\Tr( \Fr_q, D\cF_y)  
  = \varepsilon_d \varepsilon_y q^{-d} q_y^{-1} \overline{\Tr(\Fr_q,D\cF_y)},
$$
for some $d \in \Z$  (setting $\varepsilon_d = (-1)^d$). Here $D\cF$ denotes the
Verdier dual of $\cF$, and $h \mapsto \bar{h}$ is the involution of $\Z[q^{1/2},q^{-1/2}]$ determined by 
$q^{1/2} \mapsto q^{-1/2}$.

In other words, we have
$$
\Tr(\Fr_q, D\cF) = \sum_{y} \varepsilon_d \varepsilon_y q^{-d} q_y^{-1} \overline{\Tr(\Fr_q,D\cF_y)} T_y
$$
in the Iwahori-Hecke algebra.
By applying Verdier duality (i.e. by applying the Kazhdan-Lusztig involution $\overline{\cdot}$)
and using the definition $\overline{T}_y = T_{y^{-1}}^{-1}$, we get
$$
\Tr(\Fr, \cF) = \sum_{x} \left( \sum_{y \geq x} \varepsilon_d q^d \Tr
(\Fr_q, D\cF_y) \varepsilon_x R_{x,y}(q)\right) T_x.
$$
Therefore we see that property (P) implies
$$
\Tr(\Fr_q,\cF_x) = \varepsilon_d q^d \sum_{y \geq x} \varepsilon_x 
\Tr(\Fr_q,D\cF_y) R_{x,y}(q).
$$

{\em Remark.} We will see in section \ref{RPsi_satisfies_P} that
 $R\Psi$ satisfies (P). More generally, the notion of Wakimoto sheaf
allows us to define a much larger class of sheaves that satisfy (P);
cf. section \ref{wakimoto_sheaves}.

Arguing as above, it is easy to see from the definitions that property (P) is stable under Verdier duality.

\begin{lem} \label{duality_for_P}
 A sheaf $\cF \in P^{\cB}_{\sqrt q}(\Fl,\Ql)$ satisfies (P) for an integer $d$ if and only if $D\cF$ satisfies (P) for $-d$.
\end{lem}

{\em Proof.}  Use the identities
\begin{eqnarray*}
\Tr(\Fr_q, \cF_x) &=& \sum_{x \le y} \varepsilon_x \varepsilon_y q^{-1}_y R_{x,y} \overline{\Tr(\Fr_q, D\cF_y)}  \\
\overline{\Tr(\Fr_q, \cF_x)} &=& \sum_{x \le y} q_x R_{x,y} \Tr(\Fr_q, D\cF_y).
\end{eqnarray*}
\qed

\subsection{Cohomological interpretation of the multiplicities}

In this section we will give a cohomological interpretation of the
multiplicities associated to a sheaf $\cF$ which satisfies the property
(P). This can be seen as a generalization of the theorem of Kazhdan and Lusztig
which gives an interpretation of the inverse Kazhdan-Lusztig
polynomials (Theorem \ref{invKLpol}). Their theorem is the key
ingredient of the proof.

\begin{thm} \label{cohomological interpretation}
Suppose that $\cF \in P^\cB_{q}(\Fl, \Ql)$ satisfies (P).
Then
\begin{eqnarray*}
\sum_i m(\cF, w,i)q^i &=& \varepsilon_d q^d \Tr(\Fr_q, H_c^\bullet(\Fl, D\cF \tens
IC(\cB^{\overline{w}})))\\
&=& \varepsilon_d q^{d- \ell(w)} \Tr(\Fr_q, H^\bullet_c(\Fl, \cF \tens IC(\cB^{\overline{w}}))^\vee).
\end{eqnarray*}
\end{thm}
Here, $V^\vee = {\rm Hom}(V, \Ql)$ for a $\Ql$-space $V$.

{\em Proof.}
Using the definition of the multiplicities and the fact that $\cF$ satisfies
(P), we have
\begin{eqnarray*}
\sum_i m(\cF, w,i)q^i 
&=& \sum_{w \leq x} \varepsilon_x \Tr(\Fr_q,\cF_x) Q_{w,x} \\ 
&=& \sum_{w \leq x} \sum_{x \leq y} \varepsilon_d q^d \Tr(\Fr_q, (D\cF)_y)
Q_{w,x}R_{x,y}.
\end{eqnarray*}

Rearranging this and using the formula 
$\sum_{w \leq x \leq y} Q_{w,x}R_{x,y} = q_y q_w^{-1} Q_{w,y}(q^{-1})$, we get
\begin{eqnarray*}
\sum_i m(\cF, w,i)q^i &=& \varepsilon_d q^d \sum_{w \leq y} \Tr(\Fr_q,(D\cF)_y) \sum_{w \leq x \leq y} Q_{w,x}R_{x,y} \\
&=& \varepsilon_d q^d \sum_{w \leq y} (q_y q^{-1}_w Q_{w,y}(q^{-1})) \Tr(\Fr_q,(D\cF)_y).
\end{eqnarray*}

By using the Poincar\'e dual of Theorem \ref{invKLpol}, we see that 
$$
q_y q^{-1}_w Q_{w,y}(q^{-1}) = \Tr(\Fr_q, IH^\bullet_c(\cB_y \cap \cB^{\overline{w}})).
$$

Thus we get 
\begin{eqnarray*}
&& \sum_i m(\cF, w,i)q^i \\ 
&=& \varepsilon_d q^d \sum_{w \leq y} \Tr(\Fr_q, IH^\bullet_c(\cB_y \cap
\cB^{\overline{w}}))\Tr(\Fr_q,(D\cF)_y)\\
&=& \varepsilon_d q^d \sum_{w \leq y} \Tr(\Fr_q, H^\bullet_c(\cB_y \cap
\cB^{\overline{w}}, IC(\cB_y \cap \cB^{\overline{w}}) \tens D\cF|_{\cB_y \cap \cB^{\overline{w}}})) \\
&=& \varepsilon_d q^d  \Tr(\Fr_q, H^\bullet_c(\Fl, IC(\cB^{\overline{w}}) \tens D\cF))),
\end{eqnarray*}
taking Proposition \ref{restriction} into account.  This proves the first equality of the theorem.  

The second equality can be proved by applying the first equality to the sheaf $D\cF$ in place of $\cF$ (cf. Lemma \ref{duality_for_P}), and by using Lemma \ref{proto_palindromic}.
\qed

{\em Remark.} It is clear from the proof of the theorem that
we could replace the condition that $\cF$ satisfies (P) by the
following: Let $\tilde{\cF}$ be a $\cB$-equivariant sheaf such that for all $y$
$$ \Tr(\Fr_q, D\cF_y) = \varepsilon_y q_y^{-1} \overline{\Tr(\Fr_q,
  \tilde{\cF}_y)}. $$
Then
$$
\sum_i m(\cF, w,i)q^i = \Tr(\Fr_q, H_c^\bullet(\Fl, \tilde{\cF} \tens
IC(\cB^{\overline{w}}))).
$$

\subsection{The Kazhdan-Lusztig theorem as a consequence of the main theorem}

Fix $y \in \widetilde{W}$, and consider the sheaf $\cF = j_{y!}\Ql[\ell(y)]$.  This belongs to $P^{\cB}_q(\Fl,\Ql)$, and satisfies property (P) with $d = \ell(y)$.  We may therefore apply the second equality of Theorem \ref{cohomological interpretation} to this sheaf.  Using the projection formula for $j_{y!}$, we have
\begin{eqnarray*}
\varepsilon_y q_y q^{-1}_w \Tr(\Fr_q, H^\bullet_c(\Fl, \cF \tens IC(\cB^{\overline{w}}))^\vee) 
&=& q_y q_w^{-1} \Tr(\Fr_q, H^\bullet_c(\Fl, j_{y!}(j_y^*IC(\cB^{\overline{w}})))^\vee) \\
&=& q_y q_w^{-1} \Tr(\Fr_q, H^\bullet_c(\cB^{\bar{w}}_y, IC(\cB^{\overline{w}}_y))^\vee) \\
&=& q_y q_w^{-1} \Tr(\Fr^{-1}_q, H^\bullet_c(\cB^{\bar{w}}_y, IC(\cB^{\overline{w}}_y))) \\
&=& \Tr(\Fr_q, IH^\bullet(\cB_y^{\overline{w}})),\\
\end{eqnarray*}
where the last equality follows from Poincar\'{e} duality.
Since $\sum m(\cF, w, i) q^i = Q_{w,y}(q)$, this is just one version of the theorem of
Kazhdan and Lusztig (Theorem \ref{invKLpol}).  (Taking into account the very-purity of $\cB^{\overline{w}}_y$, one can recover the full statement of that theorem.)

\subsection{Consequences}

\subsubsection{$w \in (\widetilde{W}/W)_{min}$}

Denote by $\pi : \Fl \lto \Grass$ the projection. The map $\pi$ is smooth and proper as a morphism of ind-schemes.
We denote the projection $\widetilde{W} \lto \widetilde{W}/W$ by $\pi$, too.

If $w \in (\widetilde{W}/W)_{min}$, we have $\cB^{\overline{w}} =
\pi^{-1}(Q^{\overline{\pi(w)}})$, and so $IC(\cB^{\bar{w}}) = \pi^*IC({\mathcal Q}^{\overline{\pi(w)}})$.  
Here ${\mathcal Q}^{\overline{\pi(w)}}$ is the closure of the $\cN^-$-orbit of $\pi(w)$ in the affine Grassmannian for $G$.   

Thus in this case we can simplify the formula in the theorem
by considering the push-forward to $\Grass$ and applying the projection
formula. We get 
\begin{eqnarray*}
\sum_i m(\cF, w,i)q^i & = & \varepsilon_d q^d \Tr(\Fr_q, H_c^\bullet(\Fl, D\cF \tens
IC(\cB^{\overline{w}}))) \\
& = & \varepsilon_d q^d \Tr(\Fr_q, H_c^\bullet(\Grass, (\pi_\ast D\cF) \tens
IC(\cQ^{\overline{\pi(w)}}))). 
\end{eqnarray*}

We remark that this formula simplifies further in the case $\cF = R\Psi$; see section 6.3.

\section{The inertia action on $R\Psi$}
\label{inertia_action}

Elaborating an argument due to Gaitsgory \cite{Gaitsgory}, in this section we prove that $R\Psi$ carries a purely unipotent action of the inertia group.  This has two important consequences.  First, it implies that the Jordan-H\"{o}lder series for $R\Psi$ in the category $P^{\mathcal B}(\Fl \times_s \eta, 
\overline{\mathbb Q}_\ell)$ is determined by that of its image in the category 
$P^{\mathcal B}_{\rm Weil}(\Fl,\overline{\mathbb Q}_\ell)$ (and that the choice of lift of Frobenius is immaterial for the latter).  Secondly, it shows that the semi-simple trace of Frobenius and the usual trace of Frobenius agree on these nearby cycles, so that the Kottwitz conjecture can effectively be used to derive information about the Jordan-H\"{o}lder series in this latter category.  (See sections 6,7 and 9.)

\subsection{Unipotent/non-unipotent decompositions}

Let $(S,\eta,s)$ be a Henselian trait, and let $(\bar{S},\bar{\eta},\bar{s})$ be the trait constructed by letting $\bar{\eta}$ be a separable closure of $\eta$, $\bar{S}$ the normalization of $S$ in $\bar{\eta}$, and $\bar{s}$ the corresponding closed point.  We assume the residue field $k(s) = \F_q$, $q$ a power of a prime $p$.    
Write $\Gamma = {\rm Gal}(\bar{\eta}/\eta)$, $\Gamma^s = {\rm Gal}(\bar{s}/s)$, 
and let $\Gamma_0$ denote the inertia subgroup, i.e., the kernel of the canonical 
homomorphism $\Gamma \rightarrow \Gamma^s$.
For a fixed prime $l \neq p$, there is a canonical surjective homomorphism
$$
t_l: \Gamma_0 \rightarrow \Z_l(1).
$$

Let $\Q_l$-Mod denote the category of finite dimensional $\Q_l$-vector spaces.  
We consider the abelian category ${\rm Rep}(\Gamma,\Q_l)$ of continuous finite-dimensional $l$-adic representations $(\rho,V)$
$$
\rho: \Gamma \rightarrow {\rm GL}(V).
$$
By a theorem of Grothendieck (cf. \cite{Serre-Tate}), the restricted representation $\rho(\Gamma_0)$ is quasi-unipotent i.e., there exists a finite-index subgroup $\Gamma_1$ of $\Gamma_0$ which acts unipotently on $V$.  There exists then a unique nilpotent morphism, the {\em logarithm of the unipotent part of} $\rho$
$$
N : V(1) \rightarrow V
$$ 
characterized by the following property: for all $g \in \Gamma_1$, we have
$$
\rho(g) = {\rm exp}(Nt_l(g)).
$$

We make a choice of geometric Frobenius $\Phi$ lifting the automorphism $a \mapsto a^{1/q}$ with respect to 
the homomorphism ${\rm Gal}(\bar{\eta}/\eta) \rightarrow {\rm Gal}(\bar{s}/s)$.  Then $\rho$ uniquely 
determines a representation $\tilde{\rho}:\langle \Phi \rangle \ltimes \Gamma_0 \rightarrow 
{\rm GL}(V)$, trivial on some finite-index subgroup $\Gamma_1$ of $\Gamma_0$, by the formula   
$$
\rho(\Phi^n \sigma) = \tilde{\rho}(\Phi^n \sigma){\rm exp}(Nt_l(\sigma))  
$$
for $n \in \Z$ and $\sigma \in \Gamma_0$ (cf. \cite{Deligne-Antwerp}); moreover we have the relation 
$\tilde{\rho}(\Phi^n \sigma)N \tilde{\rho}(\Phi^n \sigma)^{-1} = q^{-n} N$.

We can define an idempotent endomorphism $\gamma = \gamma_V \in {\rm End}_{{\mathbb Q}_\ell}(V)$ by 
$$
\gamma = |\Gamma_0/\Gamma_1|^{-1}\sum_{g \in \Gamma_0/\Gamma_1} \tilde{\rho}(g).
$$ 
Clearly this is independent of the choice of coset representatives for $\Gamma_0/\Gamma_1$ as well as the choice of finite-index subgroup $\Gamma_1 \subset \Gamma_0$ on which $\tilde{\rho}$ is trivial.  Using this 
it is easy to see that $\gamma$ commutes with 
$\tilde{\rho}(\langle \Phi \rangle \ltimes \Gamma_0)$ and $N$, hence also with 
$\rho(\Gamma)$.  We obviously have
$$
\gamma^2 = \gamma.
$$

We note that if $f: V \rightarrow W$ is a morphism in ${\rm Rep}(\Gamma,\Q_l)$, then $f \circ \gamma_V = \gamma_W \circ f$.

\medskip

We define the {\em unipotent} and {\em non-unipotent} 
subrepresentations of $V$ as follows: 
$V^{un}$ (resp. $V^{non-un}$) is the largest $\Gamma$-subrepresentation such that, 
as a representation of $\Gamma_0$, all of its irreducible subquotients are trivial 
(resp. non-trivial).  It is easy to see that these exist, and that $V^{un} \cap V^{non-un} = (0)$.  

We denote by ${\rm Rep}(\Gamma,\Q_l)^{un}$ (resp. ${\rm Rep}(\Gamma,\Q_l)^{non-un}$) the maximal full subcategory of ${\rm Rep}(\Gamma,\Q_l)$ whose objects are all unipotent (resp. non-unipotent).

\begin{lem}\label{vector-decomp}
\begin{enumerate}
\item[(1)]  For every $l$-adic representation $(\rho,V)$ of $\Gamma$, there is a canonical 
decomposition 
$$
V = V^{un} \oplus V^{non-un}.
$$  
Moreover, any morphism $f:V \rightarrow W$ in ${\rm Rep}(\Gamma,\Q_l)$ decomposes canonically as $f = f^{un} \oplus f^{non-un} : V^{un} \oplus V^{non-un} \rightarrow W^{un} \oplus W^{non-un}$.  
This gives rise to a decomposition of abelian categories 
$$
{\rm Rep}(\Gamma,\Q_l) = {\rm Rep}(\Gamma,\Q_l)^{un} \oplus {\rm Rep}(\Gamma,\Q_l)^{non-un}.
$$ 
\item[(2)]  The endomorphism $\gamma$ is $\rho(\Gamma)$-equivariant.  Furthermore, we have $V^{un} = {\rm ker}(1 - \gamma) = {\rm im}(\gamma)$ and 
$V^{non-un} = {\rm ker}(\gamma) = {\rm im}( 1- \gamma)$.
\item[(3)] $V$ is unipotent (resp. non-unipotent) if and only if $1 - \gamma$ (resp. $\gamma$) acts by zero on $V$.
\end{enumerate} 
\end{lem}

{\em Proof.}
(1).  We only need to show that $V$ has {\em some} decomposition $V = V_1 \oplus V_2$, where $V_1$ is unipotent and $V_2$ is non-unipotent.
    
Since the restriction of $\tilde{\rho}$ to $\Gamma_0$ is semi-simple, we may  consider the decomposition of $V$ according to its isotypical components:
$$
V = \bigoplus_\tau V_\tau.
$$

For any fixed $\sigma$, the automorphism ${\rm exp}(N t_l(\sigma))$ commutes with every $\tilde{\rho}(\sigma')$, 
$\sigma' \in \Gamma_0$, hence preserves each component $V_\tau$.  Thus $\rho(\sigma)$ preserves each component $V_\tau$, for $\sigma \in \Gamma_0$.  On the other hand, it is easily seen that $\rho(\Phi)$ preserves the isotypical component $V_{\mathbb I}$ corresponding to the trivial representation, while permuting the components $V_\tau$ with $\tau \neq {\mathbb I}$.  Since the Weil group $\langle \Phi \rangle \ltimes \Gamma_0$ is dense in $\Gamma$, it follows that $V_1 = V_{\mathbb I}$ and $V_2 = \bigoplus_{\tau \neq {\mathbb I}} V_\tau$ have the required properties.  

\smallskip

(2).   We have already noted that $\gamma$ is $\rho(\Gamma)$-equivariant.  The exact sequence
$$ \xymatrix{
0 \ar[r] & V^{\tilde{\rho}(\Gamma_0)} \ar[r]^i & V \ar[r]^p &
V/V^{\tilde{\rho}(\Gamma_0)} \ar[r] & 0
} $$
splits in the category of $\rho(\Gamma)$-modules; in fact we have $\gamma \circ i = {\rm id}$ and $p \circ (1 - \gamma) = {\rm id}$.  
This shows that 
$$
V^{un} = V^{\tilde{\rho}(\Gamma_0)} = {\rm im}(\gamma) 
$$
$$
V^{non-un} \cong V/V^{\tilde{\rho}(\Gamma_0)} \cong {\rm ker}(\gamma).
$$

\smallskip

(3).  This follows from (2).
\qed

\bigskip

Let $\pi_1$ denote a pro-finite group (e.g. the algebraic fundamental group of a locally Noetherian scheme at a geometric point).  Let ${\rm Rep}(\Gamma, \pi_1, \Q_l)$ denote the category 
whose objects $V$ are $\Q_l$-modules equipped with continuous commuting actions of $\pi_1$ and 
$\Gamma$ and whose morphisms are compatible with the actions.  
Then Lemma \ref{vector-decomp} implies that this category decomposes into unipotent and 
non-unipotent subcategories.

\begin{lem}\label{fundgrp-decomp}
There is a canonical decomposition of abelian categories
$$
{\rm Rep}(\Gamma,\pi_1,\Q_l) = {\rm Rep}(\Gamma,\pi_1, \Q_l)^{un} \oplus 
{\rm Rep}(\Gamma, \pi_1, \Q_l)^{non-un}.
$$
Similarly, there is a decomposition
$$
{\rm Rep}(\Gamma_0, \pi_1,\Q_l) = {\rm Rep}(\Gamma_0, \pi_1,\Q_l)^{un} \oplus
{\rm Rep}(\Gamma_0,\pi_1, \Q_l)^{non-un},
$$
and the restriction functor ${\rm Rep}(\Gamma,\pi_1, \Q_l) \rightarrow 
{\rm Rep}(\Gamma_0,\pi_1,\Q_l)$ is compatible with the decompositions.
\qed
\end{lem}

Let $X$ denote a locally Noetherian scheme.  Let $LS(X,\Q_l)$ denote the category of locally constant and constructible $\Q_l$-sheaves on $X$, i.e., the $l$-{\em adic local systems on} $X$ (see \cite{SGA 5}, Ex. VI $\S 1.4$).  An object 
in this category is a surjective projective system 
${\mathcal F} = ({\mathcal F}_n)_{n \in \N}$ where ${\mathcal F}_n$ is a finite and locally 
constant $\Z/l^n\Z$-sheaf on $X_{\text{\'{e}t}}$.  If $X$ is connected, then by loc.cit. $\S1.4.2$, there is an equivalence of 
categories
$$
LS(X,\Q_l) \cong {\rm Rep}(\pi_1(X,a), \Q_l),
$$ 
given by the functor $F_a : ({\mathcal F}_n)_n \mapsto 
\varprojlim({\mathcal F}_n)_a \tens_{\Z_\ell} \Q_\ell$, where 
  ${\mathcal G}_a$ is the fiber of a sheaf ${\mathcal G}$ at a geometric point $a$ of $X$.  (The choice of $a$ is essentially irrelevant, since $X$ is connected.)  

Now suppose $X$ is a scheme of finite-type over $s = \Spec(\F_q)$, and
denote the geometric special fiber by 
$X_{\bar{s}} = X \otimes_{\F_q} \bar{\F}_q$.  If $X_{\bar{s}}$ is connected, we 
denote its fundamental group by $\bar{\pi}_1 = \pi_1(X_{\bar{s}})$.

We will define several categories attached to $X/s$.

\begin{itemize}
\item[-] $LS(X \times_s \eta, \Q_l)$ : objects are objects in $LS(X_{\bar{s}}, \Q_l)$ equipped with a continuous action of $\Gamma$ compatible with the action of $\Gamma$ on $X_{\bar{s}}$ (the quotient $\Gamma^s$ of $\Gamma$ acts on 
$X_{\bar{s}}$ by transport of structure).  In other words, $({\mathcal F}_n)_n$ is an object of this category if for each $g \in \Gamma$ there is an isomorphism of projective systems
$$
\phi_g : \bar{g}_*{\mathcal F}_n \cong {\mathcal F}_n
$$
such that $\phi_{gh} = \phi_g \circ \bar{g}_*(\phi_h)$.  Here $\bar{g}$ denotes the image of $g$ under the projection $\Gamma \rightarrow \Gamma^s$.  
(Note: The ``continuity'' condition is that for each quasi-compact \'{e}tale morphism $U \rightarrow X_s$ and corresponding base change 
$\bar{U} \rightarrow X_{\bar{s}}$, the action of $\Gamma$ on the (finite) set 
${\mathcal F}_n(\bar{U})$ is continous for each $n$ (cf. \cite{SGA 7}, Ex. XIII).)     

Let $LS(X \times_s \eta^{nr},\Q_l)$ denote the category consisting  of local systems with an 
action of $\Gamma_0$ rather than a compatible action of the group $\Gamma$.  
If $X_{\bar{s}}$ is connected, then taking stalks 
${\mathcal G} \mapsto {\mathcal G}_a$ gives an equivalence of categories
$$
LS(X \times_s \eta^{nr},\Q_l) \cong {\rm Rep}(\Gamma_0, \bar{\pi}_1, \Q_l).
$$ 

(Remark: This means that to define the notions of unipotent and non-unipotent for objects on the left hand side, we should consider their stalks, where such notions are already defined; cf. Lemma \ref{fundgrp-decomp}.)

\item[-] $Sh_c(X \times_s \eta,\Q_l)$ :  Let $Sh_c(X_{\bar{s}}, \Q_l)$ denote the abelian category consisting of $\Q_l$-sheaves 
$({\mathcal F}_n)_n$ on $X_{\bar{s}}$ where each ${\mathcal F}_n$ is finite and constructible.  
   
We define $Sh_c(X \times_s \eta, \Q_l)$: the objects are objects of $Sh_c(X_{\bar{s}}, \Q_l)$ equipped  with a continuous action of $\Gamma$, compatible with the action of $\Gamma$ on $X_{\bar{s}}$, in the sense defined above.  The morphisms are compatible with the action of $\Gamma$.  Continuity has the same meaning as in $LS(X \times_s \eta,\Q_l)$ above.  

The category $Sh_c(X \times_s \eta^{nr}, \Q_l)$ is defined analogously to $LS(X \times_s \eta^{nr}, \Q_l)$.

\item[-] $P(X \times_s \eta, \Q_l)$ : the objects are middle perverse sheaves on $X_{\bar{s}}$ equipped with a continuous action of $\Gamma$ compatible with its action on $X_{\bar{s}}$.  (Here continuity is tested on cohomology sheaves.)  Morphisms are compatible with the action of $\Gamma$.  This category is naturally a full subcategory of the category $D^b_c(X \times_s \eta, \Q_l)$ defined below.  

The category $P(X \times_s \eta^{nr}, \Q_l)$ is defined analogously to $LS(X \times_s \eta^{nr}, \Q_l)$.

\item[-] $D^b_c(X \times_s \eta, \Q_l)$ : the objects are objects of the ``derived'' category $D^b_c(X_{\bar{s}}, \Q_l)$ equipped with a continuous action of $\Gamma$ compatible with its action on $X_{\bar{s}}$.  
(Continuity is tested on the cohomology sheaves.)  The morphisms are compatible with the action of $\Gamma$.

The category $D^b_c(X \times_s \eta^{nr}, \Q_l)$ is analogous to 
$LS(X \times_s \eta^{nr}, \Q_l)$.

\end{itemize}

\bigskip

The following lemma extends certain properties of ${\rm Rep}(\Gamma_0,\Q_l)$ to each of the analogous categories above.

\begin{lem} \label{general-quasi-unipotence}
Let $M$ be an object of one of the categories $LS(X \times_s \eta^{nr}, \Q_l)$, 
$Sh_c(X \times_s \eta^{nr}, \Q_l)$, or $D^b_c(X \times_s \eta^{nr}, 
Q_l)$.  Then
\begin{itemize}
\item[(1)]  The action of $\rho(\Gamma_0)$ on $M$ is quasi-unipotent: there exists a finite-index subgroup $\Gamma_1 \subset \Gamma_0$ such that $(\rho(g) - 1)$ acts nilpotently on $M$, for every $g \in \Gamma_1$.
\item[(2)]  There is a unique nilpotent morphism
$N: M(1) \rightarrow M$ characterized by the equality
$$
\rho(g) = \exp(N t_l(g))
$$
for every $g \in \Gamma_1$.
\item[(3)]  There is a unique homomorphism $\tilde{\rho} : \Gamma_0 \rightarrow {\rm Aut}(M)$ which satisfies
\begin{enumerate}
\item[(a)] $\rho(\sigma) = \tilde{\rho}(\sigma) {\rm exp}(N t_l(\sigma))$, for $\sigma \in \Gamma_0$,
\item[(b)] $\tilde{\rho}$ is trivial on a finite-index subgroup of $\Gamma_0$,
\item[(c)] $\tilde{\rho}(\sigma) N \tilde{\rho}(\sigma)^{-1} = N$, for $\sigma \in \Gamma_0$.
\end{enumerate}

\end{itemize}
\end{lem}

This is closely related to the next result, whose proof works along the same lines.  To state it we need a suitable notion of ``unipotent'' and ``non-unipotent'' object for the categories above.

\begin{Def} \label{general-def-of-unipotent}
\begin{enumerate}
\item[(1)]  Let $M$ be an object of $LS(X \times_s \eta, \Q_l)$ or $Sh_c(X \times_s \eta,\Q_l)$.  
We say $M$ is {\em unipotent} (resp. {\em non-unipotent}) if for every geometric point $\bar{x} \in X_{s}(\bar{\F}_q)$, the stalk $M_{\bar{x}}$ is a unipotent (resp. non-unipotent) object of ${\rm Rep}(\Gamma_0,\Q_l)$.
\item[(2)]  Let $M$ be an object of $P(X \times_s \eta,\Q_l)$ or $D^b_c(X \times_s \eta,\Q_l)$.  We say $M$ is {\em unipotent} (resp. {\em non-unipotent}) if every cohomology sheaf ${\mathcal H}^iM$ is a unipotent (resp. non-unipotent) object of $Sh_c(X \times_s \eta, \Q_l)$.
\end{enumerate}
\end{Def}

Remark:  In the same way we can define the notions of unipotent and non-unipotent for objects of $LS(X \times_s \eta^{nr},\Q_l)$, $Sh_c(X \times_s \eta^{nr},\Q_l)$, etc.  It is clear that $M \in LS(X \times_s \eta,\Q_l)$ etc. is unipotent (resp. non-unipotent) if and only if it has unipotent (resp. non-unipotent) image under the forgetful functor 
$LS(X \times_s \eta,\Q_l) \rightarrow LS(X \times_s \eta^{nr},\Q_l)$ etc.

\begin{lem} \label{general-decomp}
\begin{enumerate}

\item[(1)]  
Let ${\mathcal C}$ be one of the categories $LS(X \times_s \eta^{nr},\Q_l)$, $Sh_c(X \times_s \eta^{nr},\Q_l)$, or $D^b_c(X \times_s \eta^{nr},\Q_l)$.  Let $M$ be an object of ${\mathcal C}$.  
Define the idempotent $\gamma = \gamma_M \in {\rm End}_{{\mathcal C}}(M)$ by
$$
\gamma = |\Gamma_0/\Gamma_1|^{-1} \sum_{g \in \Gamma_0/\Gamma_1} \tilde{\rho}(g)
$$
where $\tilde{\rho}$ is as above and $\Gamma_1 \subset \Gamma_0$ is any finite-index subgroup of $\Gamma_0$ on which $\tilde{\rho}$ is trivial.
Then $M$ is unipotent (resp. non-unipotent) if and only if $1 - \gamma_M$
(resp. $\gamma_M$) acts by zero on $M$.  
\item[(2)]  If $M$ is an object of one of the abelian categories $LS(X \times_s \eta^{nr},\Q_l)$, $Sh_c(X \times_s \eta^{nr},\Q_l)$, or $P(X \times_s \eta^{nr},\Q_l)$
 then the canonical decomposition
$$
M = {\rm ker}(1 - \gamma_M) \oplus {\rm ker}(\gamma_M)
$$
is a canonical decomposition of $M$ into unipotent and non-unipotent direct factors.
The same holds for the categories $LS(X \times_s \eta, \Q_\ell)$, 
 $Sh_c(X \times_s \eta, \Q_\ell)$, and $P(X \times_s \eta, \Q_\ell)$.
\end{enumerate}
\end{lem}

{\em Proof of Lemmas \ref{general-quasi-unipotence} and \ref{general-decomp}.}  

To begin, we assume $M \in D^b_c(X \times_s \eta^{nr},\Q_l)$, the most general category in Lemma \ref{general-quasi-unipotence}.

Once it is known that $\rho(\Gamma_0)$ acts quasi-unipotently on $M$, the existence and uniqueness of $N$ follow in general from the same argument used in the category ${\rm Rep}(\Gamma_0,\Q_l)$.  Indeed, let $P_l = {\rm ker}(t_l)$ and let $\Gamma_1 \subset \Gamma_0$ denote a finite-index subgroup such that $\rho(\Gamma_1)$ acts unipotently on $M$.  One shows that $\rho(\Gamma_1 \cap P_l) = {\rm id}_M$ (this can be tested on cohomology stalks, effectively reducing to the case $M \in {\rm Rep}(\Gamma_0,\Q_l)$ where it is already known), and thence that $\rho|\Gamma_1$ factors through $t_l$; then define $N: M(1) \rightarrow M$ by $Nt_l(\sigma) = {\rm log}(\rho(\sigma))$ for all $\sigma \in \Gamma_1$.

Moreover, given $N$ we then define $\tilde{\rho}$ by the equality
$$
\tilde{\rho}(\sigma)  = \rho(\sigma) {\rm exp}(-Nt_l(\sigma))
$$
for $\sigma \in \Gamma_0$.
The identities proving that $\tilde{\rho}$ is a homomorphism and the 
equality in Lemma \ref{general-quasi-unipotence} (3(c))
reduce to the case $M \in {\rm Rep}(\Gamma_0,\Q_l)$.  
Indeed, passing to cohomology sheaves, we may assume 
$M \in Sh_c(X \times_s \eta^{nr},\Q_l)$.  Then the constructibility of $M$ 
allows us to view $M$ as the result of  
``gluing'' a local system on an open subscheme with a constructible 
sheaf supported on a proper closed subscheme (for the precise statement of gluing, see for example \cite{Milne}, II $\S 3$).   
By Noetherian induction therefore, we may assume $M \in LS(X \times_s \eta^{nr},\Q_l)$.  Passing to the finitely many connected components of $X$, we are reduced to $M \in {\rm Rep}(\Gamma_0,\bar{\pi}_1,\Q_l)$, where the desired identities are 
already known.

To prove that $\tilde{\rho}$ is trivial on a finite-index subgroup $\Gamma_1 \subset \Gamma_0$, one uses the same reduction steps, keeping in mind that for $M \in Sh_c(X \times_s \eta^{nr},\Q_l)$, there is a finite stratification $X_i$ ($i \in I$) of $X$ by locally closed subschemes $X_i$ such that $M|X_i \in LS(X_i \times_s \eta^{nr},\Q_l)$, and each $X_i$ has only finitely many connected components.

To show that $\rho(\Gamma_0)$ acts quasi-unipotently on $M$ is similar.  This completes the proof of Lemma \ref{general-quasi-unipotence}.

Now again suppose $M$ is an object of $D^b_c(X \times_s \eta^{nr}, \Q_l)$.  We wish to prove $M$ is unipotent (resp. 
non-unipotent) if and only if $1 - \gamma_M$ (resp. $\gamma_M$) acts by zero on $M$.  Passing to cohomology sheaves and using the obvious identity
$$
\gamma_{{\mathcal H}^iM} = {\mathcal H}^i(\gamma_M)
$$
we can assume $M \in Sh_c(X \times_s \eta^{nr},\Q_l)$.  Taking stalks and using the identity
$$
\gamma_{i^*_{\bar{x}}M} = i^*_{\bar{x}}(\gamma_M),
$$
where $i_{\bar{x}} : \bar{x} \rightarrow X_s$ is the embedding of a geometric point in $X_s(\bar{\F}_q)$, we can assume that $M \in {\rm Rep}(\Gamma_0,\Q_l)$, where the result follows from Lemma \ref{vector-decomp}.  This proves 
Lemma \ref {general-decomp} part (1).

Concerning part (2), 
if $M$ belongs to one of the abelian categories $LS(X \times_s \eta^{nr}, \Q_l)$, 
$Sh_c(X \times_s \eta^{nr},\Q_l)$, or $P(X \times_s \eta^{nr},\Q_l)$, then the 
idempotent endomorphism $\gamma_M$ gives rise to the required decomposition.  
This proves part (2), and thus Lemma \ref{general-decomp}.
\qed

Remark: If ${\mathcal C}$ is one of the categories in Lemma \ref{general-decomp} (2), there is a decomposition of abelian categories
$$
{\mathcal C} = {\mathcal C}^{un} \oplus {\mathcal C}^{non-un},
$$
and we denote the corresponding components of an object $M$ by $M^{un}$ and $M^{non-un}$.

This argument does not prove that the category 
$D^b_c(X \times_s \eta^{nr},\Q_l)$ decomposes in a like manner: since it is not abelian, we cannot even define the kernel and image of $\gamma_M$ for a general object $M$.

\bigskip

\subsection{Permanence properties of unipotence and non-unipotence}

Suppose $f: X \rightarrow Y$ is a finite-type morphism of 
finite-type $S$-schemes.  This gives rise to the ``derived'' functors 
 $f_!: D^b_c(X \times_s \eta^{nr},\Q_l) \rightarrow 
D^b_c(Y \times_s \eta^{nr}, \Q_l)$ and 
 $f^*: D^b_c(Y \times_s \eta^{nr},\Q_l) \rightarrow 
D^b_c(X \times_s \eta^{nr},\Q_l)$.

\begin{lem} \label{permanence}
\begin{enumerate}
\item[(1)]  Let  $M \in D^b_c(X \times_s \eta^{nr},\Q_l)$ be unipotent (resp. non-unipotent).  Then $f_!M \in D^b_c(Y \times_s \eta^{nr},\Q_l)$ is unipotent (resp. non-unipotent).
\item[(2)]  Let $M' \in D^b_c(Y \times_s \eta^{nr},\Q_l)$, and suppose $f_{\bar{s}}: X_{\bar{s}} \rightarrow Y_{\bar{s}}$ is surjective.  Then $M'$ is unipotent (resp. non-unipotent) if and only if $f^*M'$ is unipotent (resp. non-unipotent).
\end{enumerate}
\end{lem}

{\em Proof.}
\noindent (1).  By the local-global spectral sequence
$$
{\mathcal H}^p(f_!{\mathcal H}^qM) \Rightarrow {\mathcal H}^{p+q}f_!M
$$
we are reduced to the case $M \in Sh_c(X \times_s \eta^{nr},\Q_l)$.  Taking cohomology stalks we know that $f_!M$ is unipotent if and only if for every $\bar{y} \in Y_s(\bar{\F}_q)$, and every $i$, $1 - \gamma$ acts by zero on $H^i_c(f^{-1}(\bar{y}), M)$.  But as $M$ is unipotent, $1 - \gamma$ acts by zero on $M$, hence on the cohomology group with coefficients in $M$ as well.  The non-unipotent case is similar.

\noindent (2).  This is straightforward.
\qed

\medskip

\subsection{Unipotence/non-unipotence of products and convolutions}
\label{convolve}

The following facts follow easily from Lemmas \ref{general-decomp} and \ref{permanence}.

Let ${\mathcal F}$ and ${\mathcal G}$ be  objects of $D^b_c(X \times_s
\eta,\Q_l)$ 
and suppose that ${\mathcal F}$ is unipotent.  
Then $({\mathcal F} \boxtimes {\mathcal G})^{non-un} = {\mathcal F} \boxtimes {\mathcal G}^{non-un}$.  

Now let $\mathcal B$ denote the Iwahori-group scheme acting on the affine flag
variety ${\mathcal Fl} = G(\bar{{\mathbb F}}_p((t)))/\mathcal B$.  
If the objects ${\mathcal F}$ and ${\mathcal G}$ above are $\mathcal B$-equivariant
perverse sheaves on ${\mathcal Fl}$, then we may define their twisted
product ${\mathcal F} \tilde{\boxtimes} {\mathcal G}$ and their convolution
product ${\mathcal F} \star_s {\mathcal G}$ (cf. \cite{HN1}).
Still under the assumption ${\mathcal F} = {\mathcal F}^{un}$, one can prove first $({\mathcal F} \tilde{\boxtimes} {\mathcal G})^{non-un}
= {\mathcal F} \tilde{\boxtimes} {\mathcal G}^{non-un}$, and 
then $({\mathcal F} \star_s {\mathcal G})^{non-un} = {\mathcal F} \star_s {\mathcal G}^{non-un}$.

\subsection{$R\Psi = R\Psi^{un}$}\label{nearbycyclesareunipotent}

First, we remark that in the above discussion in 5.1-5.3, we can without any problems replace
$\Q_\ell$ by a finite extension, and can then also pass to the 2-inductive
limit of the categories $D^b_c(X, E)$, $E/\Q_\ell$ finite. Thus all the
statements above hold for $\Ql$-adic sheaves as well.

We now come back to the affine flag variety. Fix a dominant coweight $\mu$
and let $R\Psi = R\Psi(IC({\mathcal Q}_\mu)[\ell(\mu)]) \in P^\cB(\Fl \times_s \eta, \Ql)$.

\begin{thm} {\rm (Gaitsgory)} The action of $\Gamma_0$ on the sheaf $R\Psi$
 of nearby cycles is purely unipotent. 
\end{thm}

{\em Proof.}
 Because $R\Psi = R\Psi^{un} \oplus R\Psi^{non-un}$ it suffices to prove that $R\Psi^{non-un} = 0$.

Note that $R\Psi^{non-un}$ is a $\mathcal B$-equivariant object of $P({\mathcal Fl} \times_s \eta, \Q_l)$
(for equivariance, use Lemma \ref{permanence}).  By using \ref{convolve} we
see that, as for $R\Psi$ itself, $R\Psi^{non-un}$ is central with respect to
convolution of Iwahori-equivariant perverse sheaves on the affine flag
variety.  To show it is zero, it suffices
to show that its semi-simple trace function $\tau_{R\Psi^{non-un}} : x
\mapsto {\rm Tr}^{ss}({\rm Fr}_q, R\Psi^{non-un}_x)$, an element in 
$Z({\mathcal H}_I)$, vanishes.  

Now consider the projection $\pi: \Fl \rightarrow \Grass$.  
Because the semi-simple trace provides a good sheaf-function
dictionary in the sense of Grothendieck (cf. \cite{HN1}, Proposition 10), the
semi-simple trace function $\Tr^{ss}(\Fr_q, \pi_{*}(R\Psi^{non-un}))$ 
corresponds to summing
$\Tr^{ss}(\Fr_q, R\Psi^{non-un})$ 
along the fibers of $\pi$, i. e., we have the following equality of 
functions in the
spherical Hecke algebra $\cH_K$: 
$$
\Tr^{ss}(\Fr_q, \pi_{*}(R\Psi^{non-un})) = \Tr^{ss}(\Fr_q,
R\Psi^{non-un}) \star {\mathbb I}_K.
$$ 
By a theorem of Bernstein (cf. \cite{L1}), 
 $- \star {\mathbb I}_K : Z({\mathcal H}_I) \tilde{\rightarrow} 
 {\mathcal H}_K$ is an isomorphism, 
so it suffices to prove $\pi_{*}(R\Psi^{non-un}) = 0$. 
But by Lemma \ref{permanence} we have
$$
\pi_{*}(R\Psi^{non-un}) = (\pi_{*}R\Psi)^{non-un} = 
  ({\mathcal A}_{\mu,s})^{non-un} = 0.
$$
Here ${\mathcal A}_{\mu}$ is the ``constant'' intersection complex $IC({\mathcal Q}_\mu)[\ell(\mu)]$ on the ``constant family'' ${\rm Grass}$, hence it has trivial inertia action.  We are done.

\subsection{Semi-simple trace vs. usual trace}

\begin{stz}
The trace and the semi-simple trace on $R\Psi$ coincide:
$$
{\rm Tr}^{ss}({\rm Fr}_q, R\Psi) = {\rm Tr}({\rm Fr}_q, R\Psi^{ss}).
$$
\end{stz}

Here, on the left hand side $R\Psi$ is regarded as an object of the category $P^{\mathcal B}(\Fl \times_s \eta, \overline{\mathbb Q}_\ell)$, and the semi-simple trace is defined as in \cite{HN1}.  On the right hand side, $\Fr_q$ represents a choice of lift of geometric Frobenius, and $R\Psi^{ss}$ denotes the semi-simplification of $R\Psi$ when regarded as an object of $P^{\mathcal B}_{\rm Weil}(\Fl, \Ql)$ (via that choice of lift).

{\em Proof.} It suffices to find a Galois-invariant filtration on $R\Psi$ such that $\Gamma_0$ acts trivially on the 
subquotients.  Since $\Gamma_0$ acts unipotently on $R\Psi$ by \ref{nearbycyclesareunipotent}, such a filtration is given by the kernels of powers of the operator $N$:
$$
F^iR\Psi = {\rm ker}N^i.
$$

\subsection{A different proof}

Here is a sketch of a different argument for the result
above: after fixing a Weil sheaf structure on $R\Psi$, the usual trace of
Frobenius satisfies the same sheaf-function dictionary as the semi-simple
trace: the same proof as for semi-simple trace shows that the usual trace yields a central function in the Iwahori-Hecke algebra, characterized by 
$$
\Tr(\Fr_q, R\Psi) \star \mathbb I_K = \Tr(\Fr_q, IC({\mathcal Q}_\mu)[\ell(\mu)]).
$$
It follows from this that it coincides with the
semi-simple trace function. This in turn implies that the action of inertia
on $R\Psi$ must be purely unipotent. To prove the last statement, it is 
enough to show that $R\Psi$ admits a filtration such that the inertia group
$\Gamma_0$ acts trivially on the graded pieces. Consider the filtration 
$$ F^i R\Psi = \ker N^i,$$
where $N$ is the logarithm of the unipotent part of the local monodromy, 
as before. The inertia group acts on the $F^i/F^{i-1}$ through a finite
quotient, so the definition of semi-simple trace and the above remarks yields
$$ \Tr(\Fr_q, (\bigoplus F^i/F^{i-1})^{\Gamma_0}) =
   \Tr^{ss}(\Fr_q, R\Psi) = \Tr(\Fr_q, R\Psi) = 
\Tr(\Fr_q, \bigoplus F^i/F^{i-1}). $$
Thus the function of the (necessarily perverse) quotient 
 $(\bigoplus F^i/F^{i-1})/ (\bigoplus F^i/F^{i-1})^{\Gamma_0}$
is zero, so we have 
 $(\bigoplus F^i/F^{i-1})^{\Gamma_0} = \bigoplus F^i/F^{i-1}$.  (We have used the fact that if ${\mathcal F} \in P(\Fl \times_s \eta, \Ql)$ carries a {\em finite} action of inertia, then ${\mathcal F}^{\Gamma_0} = {\rm im}(\gamma_{\mathcal F})$, hence is also perverse.)

\section{Applying the theorem to $R\Psi$}
\label{RPsi_satisfies_P}

Fix a dominant coweight $\mu$ and 
write $R\Psi = R\Psi( IC(\cQ_\mu)[\ell(\mu)]) \in P^\cB(\Fl \times_s \eta,\Ql)$.

\subsection{$R\Psi$ as a Weil sheaf}

Fixing a lift $\Fr_q \in \Gamma$ of the inverse of the Frobenius morphism
enables us to view $R\Psi$ as a Weil-perverse sheaf.  Since it is also a ${\mathcal B}$-equivariant perverse sheaf, it is a {\em mixed} Weil-perverse sheaf (this also follows from the Appendix, Theorem \ref{variant}). 

The Kottwitz conjecture holds for $R\Psi$ (see section 2.7), so
$$
\Tr(\Fr_q, R\Psi) = \varepsilon_\mu q^{1/2}_\mu \sum_{\lambda \leq 
        \mu}m_\mu(\lambda) z_\lambda.
$$
In particular, $\Tr(\Fr_q, R\Psi_x) \in \Z[q,q^{-1}]$, so we can consider $R\Psi$ 
as an object of $P^\cB_q(\Fl, \Ql)$.

We will show that the nearby cycles sheaf
satisfies the property (P), and hence we can apply Theorem 
\ref{cohomological interpretation}, and get some information on the
Jordan-H\"older series of $R\Psi$ as a Weil sheaf.

But furthermore, we know that the inertia group $\Gamma_0$ acts purely
unipotently on $R\Psi$, so in this case, we get a stronger independence
statement with respect to the choice of the lift $\Fr_q$. Namely, we can
regard the sheaves $IC_w(-i)$, with trivial $\Gamma_0$-action, 
as objects of $P(\Fl \times_s \eta, \Ql)$, and then the Jordan-H\"older
series of $R\Psi \in P(\Fl \times_s \eta, \Ql)$ is 'the same' as that of
$R\Psi$ as a Weil-perverse sheaf.

\subsection{$R\Psi$ and property (P)}

To show that $R\Psi$ satisfies the property (P), we need to prove
\begin{equation*} 
\overline{\Tr(\Fr_q, R\Psi_x)} = 
 \varepsilon_x \varepsilon_\mu q_x q_\mu^{-1} \Tr(\Fr_q, R\Psi_x),
\end{equation*}
because $DR\Psi = R\Psi(\ell(t_\mu))$, i.e. 
 $\Tr(\Fr_q, D R\Psi) = q_\mu^{-1}\Tr(\Fr_q, R\Psi)$.

Note the following general lemma. We write $Q = q^{-1/2} - q^{1/2}$.

\begin{lem} \label{P_lemma}
Let $f \in \Z[q,q^{-1}]$, $\alpha \in \frac{1}{2}\Z$ such that $f = q^\alpha R(Q)$
for some $R \in \Z[Q]$. Then all $Q$-powers occurring in $R$ have the 
parity of $2\alpha$, and $f(q^{-1}) = (-1)^{2\alpha}
q^{-2\alpha} f(q)$. \qed
\end{lem}

Recalling that $z_\lambda = \sum_{\nu \in W(\lambda)} \Theta_\nu$ and the fact that
$\Theta_\nu = \tilde{T}_{\nu_1} \tilde{T}^{-1}_{\nu_2}$, where $\nu_1 - \nu_2 = \nu$ and $\nu_i$ is dominant, we see (using the Kottwitz conjecture) that 
$$
R\Psi(x) = \varepsilon_\mu q_\mu^{1/2} q_x^{-1/2} R(Q),
$$ 
where $R \in {\Bbb N}[Q]$. Thus we can apply the lemma with 
 $\alpha = (\ell(\mu)-\ell(x))/2$, and get that $R\Psi$ satisfies (P).

\subsection{Consequences of the theorem for $R\Psi$}

Since $R\Psi$ is self-dual up to a Tate twist, from Theorem \ref{cohomological interpretation} we get

\begin{kor} We have
$$ \sum_i m(R\Psi, w, i)q^i = \varepsilon_\mu
   \Tr(\Fr_q, H^\bullet_c(\Fl, R\Psi \tens IC(\cB^{\overline{w}}))). $$
\end{kor}

In particular, by setting $w = \tau$, the unique element of length 0 in the
support of $R\Psi$, we get by using the projection formula for $\pi: \Fl \rightarrow {\rm Grass}$ that

\begin{eqnarray*}
\sum_i m(R\Psi, \tau, i)q^i & = &
\varepsilon_\mu \Tr(\Fr_q, H^\bullet_c(\Fl, R\Psi)) \\
& = & \varepsilon_\mu \Tr(\Fr_q, H^\bullet_c(\Grass, \pi_\ast R\Psi)) \\
& = & \varepsilon_\mu \Tr(\Fr_q, H^\bullet_c(\Grass, IC(\cQ_\mu)[\ell(\mu)])) \\
& = & \sum_i \dim IH^{2i}(\cQ_\mu, \Ql) q^i. 
\end{eqnarray*}

\begin{kor} \label{remarkable_fact} The multiplicity $m(R\Psi, \tau, i)$ is the $2i$-th
  intersection Betti number of $\bar{\cQ}_\mu$.
\end{kor}

{\em Remark.}
This corollary was initially noticed as a purely empirical fact with the aid of a computer program, and it was the starting point of this work: the effort to explain it conceptually lead to the discovery of the more general Theorem \ref{cohomological interpretation}.  The corollary can also be proved, taking the Kottwitz conjecture into account, by purely combinatorial
considerations in the Iwahori-Hecke algebra $\cH$ (recall that the traces of
Frobenius and the multiplicities are related by a base change in $\cH$).
But it seems impossible to prove anything about the multiplicities for
other $w$'s in that way.

More generally, by applying section 4.6.1 to ${\mathcal F} = R\Psi$, we get 
\begin{kor} Suppose $w \in {\rm Adm}(\mu)$ belongs to $(\widetilde{W}/W)_{min}$.  Then the multiplicity function is given by 
$$
\sum_{i} m(R\Psi,w,i) q^i = \Tr(\Fr_q, H^\bullet_c(\bar{\mathcal Q}_\mu, IC(\bar{\mathcal Q}_\mu) \otimes IC({\mathcal Q}^{\overline{\pi(w)}}))).
$$
\end{kor}

\section{Wakimoto sheaves}
\label{wakimoto_sheaves}

\subsection{Definition and properties}

We have a convolution product 
$$\star \colon P^\cB_{\rm Weil}(\Fl, \Ql)
\times P^\cB_{\rm Weil}(\Fl, \Ql) \lto D^b_{c,\rm Weil}(\Fl, \Ql),
$$
see \cite{L3}, \cite{HN1}. Using the usual sheaf-function dictionary, it is
easy to see that convolution of sheaves corresponds to multiplication of
their functions in the Hecke algebra.

\begin{Def}
The Wakimoto sheaf associated to $v, w \in \tilde{W}$ is the sheaf
$M_{v,w} := j_{v!}\Ql[\ell(v)] \star j_{w\ast}\Ql[\ell(w)]$.
\end{Def}

Note that, as usual, $j_{v!}$ and $j_{w*}$ denote the {\em derived} 
functors. The function in the Iwahori-Hecke algebra $\cH$ corresponding to
 $M_{v,w}$ is $\varepsilon_v \varepsilon_w q_w T_v T_{w^{-1}}^{-1} = \varepsilon_v \varepsilon_w q_v^{1/2} 
q_w^{1/2} \tilde{T}_v \tilde{T}^{-1}_{w^{-1}}$.  As usual, we write $\tilde{T}_x = q_x^{-1/2}T_x$.

\begin{stz} {\rm Properties of Wakimoto sheaves}
\begin{enumerate}
\item[(1)] $M_{v, w}$ is perverse.
\item[(2)] $M_{v, w}$ satisfies (P), with $d = \ell(v) + \ell(w)$.
\end{enumerate}
\end{stz}

{\em Proof.} (1) This result is due to Mirkovic. 
See \cite{HP} Prop. 6.2, or \cite{AB} Thm. 5.

(2) The function corresponding to the dual of $M_{v,w}$ is $\varepsilon_v
    \varepsilon_w q^{-1}_w T_w T^{-1}_{v^{-1}} = \varepsilon_v
    \varepsilon_w q_w^{-1/2}q_{v}^{-1/2} \tilde{T}_w \tilde{T}^{-1}_{v^{-1}} $.  We have
$$
 DM_{v,w}(y) = \varepsilon_v \varepsilon_w q_w^{-1/2} q_v^{-1/2} q_y^{-1/2} R^w_{y,v}(Q)
$$
where $R^w_{y,v} \in {\mathbb N}[Q]$ is such that $Q^i$ occurs only for $i \equiv \ell(w) + \ell(v) + \ell(y) \,\, \mbox{mod \,$2$}$.  (Use, for example, Lemma \ref{P_lemma}.)  Thus,
$$ 
 \overline{DM_{v,w}(y)} = \varepsilon_v \varepsilon_w \varepsilon_y q_w q_v q_y  DM_{v,w}(y),
$$
which proves property (P) with $d = \ell(v) + \ell(w)$.
\qed

The proposition shows that there is a large class of sheaves which satisfy
the property (P).

\begin{stz}
Let $\cF \in P^\cB_{\sqrt{q}}(\Fl, \Ql)$ be a sheaf which has a filtration
whose graded pieces are of the form $M_{v_i, w_i}(-\alpha_i)$,
where $\alpha_i \in \Z$ is such that $2\alpha_i + \ell(v_i) + \ell(w_i)$ is
independent of $i$. Then $\cF$ satisfies (P), with $d = \ell(v_i) + \ell(w_i) + 2\alpha_i$.
\end{stz}

{\em Proof.} The function corresponding to $\cF$ is 
$$ \sum_i q^{\alpha_i} q_{w_i} \varepsilon_{v_i} \varepsilon_{w_i}  T_{v_i}
T_{w_i^{-1}}^{-1}, $$
and the function of the dual of $\cF$ is obtained by applying
$\overline{\cdot}$. A computation similar to the one above yields the result. 
\qed

{\em Remark.} In their recent paper \cite{AB}, Arkhipov and Bezrukavnikov
show that the sheaf of nearby cycles $R\Psi$ has a filtration 
whose graded pieces are Wakimoto sheaves.

Their theorem relies on the sheaf of nearby cycles being central, so that we
do not get a proof of (P) for $R\Psi$ which does not use the Kottwitz
conjecture. 

\subsection{The elements in the Hecke algebra corresponding to Wakimoto sheaves}

Fix $v,w \in W_{\rm aff}$, and suppose $w = s_1 \cdots s_r$ is reduced.  There are unique polynomials 
$R^v_{x,w}(Q)$ in the parameter $Q := q^{-1/2} - q^{1/2}$ and with non-negative integral coefficients, such that
$$
\tilde{T}_v \tilde{T}_{w^{-1}}^{-1} = \sum_{x \leq v|w} R^v_{x,w}(Q) \tilde{T}_x.
$$
Here the notation $x \leq v|w$ means by definition that $x$ is the terminal element of a 
$v$-distinguished subexpression with respect to $w$ (see section 5 of \cite{H1}). Recall that an $r+1$-tuple $[\sigma_0, \dots, \sigma_r]$ is $v$-distinguished if $\sigma_0 = v$, and if for each $j \geq 1$, 
\begin{enumerate}
\item $\sigma_j \in \{ \sigma_{j-1}, \sigma_{j-1}s_j \}$; and
\item $\sigma_j = \sigma_{j-1}$ only if $\sigma_{j-1} < \sigma_{j-1}s_j$.
\end{enumerate}

We are going to prove a formula for $R^v_{x,w}(Q)$.
We need first some more definitions.

For a $v$-distinguished subexpression $\underline{\sigma} = [\sigma_0, \dots, \sigma_r]$, we set $\pi(\underline{\sigma}) = \sigma_r$.  Also define
$$
{\mathcal D}(x) = \lbrace \underline{\sigma} \, | \, \pi(\underline{\sigma}) = x \rbrace.
$$
Further define integers $n(\underline{\sigma}) = | \lbrace j \, | \, \sigma_{j-1} = \sigma_j \rbrace |$, and $m(\underline{\sigma}) = | \lbrace j \, | \, \sigma_{j-1} > \sigma_j \rbrace |$.  We have the following easy lemma.

\begin{lem}
If $\underline{\sigma}$ is $v$-distinguished with respect to $w$, then
$$
\ell(\pi(\underline{\sigma})) = \ell(v) + \ell(w) - n(\underline{\sigma}) - 2m(\underline{\sigma}).
$$
\end{lem}

\medskip

Now write $w' = s_1 \cdots s_{r-1}$, so that $w' = ws_r < w$.  The equality
$$
\sum_{y \leq v|w'} R^v_{y,w'}(Q) \tilde{T}_y(\tilde{T}_{s_r} + Q) =
\sum_{x \leq v|w} R^v_{x,w}(Q) \tilde{T}_x
$$
give rise to recurrence relations for the polynomials $R^v_{x,w}$.  Using these, we can prove the following proposition by induction on $\ell(w)$.

\begin{stz}
If $x \leq v|w$, then 
$$
R^v_{x,w}(Q) = \sum_{\underline{\sigma} \in {\mathcal D}(x)} Q^{n(\underline{\sigma})} = \sum_{\underline{\sigma} \in {\mathcal D}(x)} Q^{\ell(v) + \ell(w) - \ell(x) -2m(\underline{\sigma})}.
$$
\end{stz}

Taking the parity of the exponents into account, this proposition gives a direct verification of the conclusion of Lemma \ref{P_lemma} for Wakimoto sheaves. 

{\em Proof.}
The proof is by induction on $\ell(w)$, and uses case-by-case analysis.  Write $w'$ as above and write $s_r = s$.  We consider the set ${\mathcal D}'(y)$ of $v$-distinguished subexpressions with respect to $w'$ which end with $y$.  Integers $n$ and $m$ can be attached to each of these in the same way as before.

Note that $x \leq v|w$ implies that $xs \leq v|w'$ or $x \leq v|w'$, or both.

\medskip

{\em Case 1}: $x \leq v|w$, $x \leq v|w'$, and $xs \nleq v|w'$.

\smallskip

\noindent (a) $x < xs$.  Then every $\underline{\sigma} \in {\mathcal D}(x)$ is of the form $[?,x,x]$. Given such, write $\underline{\sigma}' = [?,x]$ (the first part of $\underline{\sigma}$).  We see $n(\underline{\sigma}) = n(\underline{\sigma}') + 1$,  $m(\underline{\sigma})= m(\underline{\sigma}')$, and $R^v_{x,w} = QR^v_{x,w'}$.  Moreover, we have a bijection of sets
$S:{\mathcal D}(x) \tilde{\rightarrow} {\mathcal D}'(x)$ given by $\underline{\sigma} \mapsto \underline{\sigma}'$.

\smallskip

\noindent (b) $xs < x$.  This cannot happen, since the only possible 
$\underline{\sigma}'$ is of form $[?,x]$, and $\underline{\sigma}$ cannot have form $[?,x,x]$, since then we would need to have $x < xs$ (by definition of $v$-distinguished).

\medskip

{\em Case 2}: $x \leq v|w$, $x \nleq v|w'$, and $xs \leq v|w'$.  

\smallskip

\noindent (a) $x < xs$.
Then every $\underline{\sigma}$ is of form $[?,xs,x]$; set 
$\underline{\sigma}' = [?,xs]$.  We have $n(\underline{\sigma}) = 
n(\underline{\sigma}')$, 
and $R^v_{x,w}= R^v_{xs,w'}$.  Also, $S: {\mathcal D}(x) \tilde{\rightarrow} {\mathcal D}'(xs)$ by $\underline{\sigma} \mapsto \underline{\sigma}'$.

\smallskip

\noindent (b) $xs < x$.  Then $\underline{\sigma}$ is of form $[?,xs,x]$; let 
$\underline{\sigma}' = [?,xs]$.  We have $n(\underline{\sigma}) = 
n(\underline{\sigma}')$, and $R^v_{x,w} = R^v_{xs,ws}$, and $S: {\mathcal D}(x) \tilde{\rightarrow} {\mathcal D}'(xs)$ by $\underline{\sigma} \mapsto \underline{\sigma}'$.

\medskip

{\em Case 3}: $ x \leq v|w$, $x \leq v|w'$, and $xs \leq v|w'$.

\smallskip 

\noindent (a) $x < xs$.  We define 
${\mathcal A}(x) = \lbrace \underline{\sigma} \,\, \mbox{is of the form} \,\, [?,xs,x] \rbrace$, and 
${\mathcal B}(x) = \lbrace \underline{\sigma} \,\, \mbox{is of the form} \,\, [?,x,x] \rbrace$.
Note that ${\mathcal D}(x) = {\mathcal A}(x) \coprod {\mathcal B}(x)$.

For $\underline{\sigma}= [?,xs,x] \in {\mathcal A}(x)$, let 
$\underline{\sigma}' = [?,xs]$.  The map $\underline{\sigma} 
\mapsto \underline{\sigma}'$ gives ${\mathcal A}(x) \tilde{\rightarrow} {\mathcal D}'(xs)$, and $n(\underline{\sigma}) = 
n(\underline{\sigma}')$. 

For $\underline{\sigma} = [?,x,x] \in {\mathcal B}(x)$, let 
$\underline{\sigma}' = [?,x]$.  The map $\underline{\sigma}
\mapsto \underline{\sigma}'$ gives ${\mathcal B}(x) \tilde{\rightarrow} {\mathcal D}'(x)$, and 
$n(\underline{\sigma}) = n(\underline{\sigma}') + 1$.

In either case, we have $R^v_{x,w} = QR^v_{x,ws} + R^v_{xs,ws}$.

\smallskip

\noindent (b) $xs < x$.  Then $\underline{\sigma} = [?,xs,x]$, and setting 
$\underline{\sigma}' = [?,xs]$, we have $n(\underline{\sigma}) = 
n(\underline{\sigma}')$, $S: {\mathcal D}(x) \tilde{\rightarrow} {\mathcal D}'(xs)$.  We have
$R^v_{x.w} = R^v_{xs,ws}$.

\medskip

The proposition now follows easily by using the induction hypothesis and by taking into account the various cases.
\qed

\medskip

In general it is difficult to say anything really concrete concerning the polynomial $R^{v}_{x,w}(Q)$: even its degree is mysterious.  However, suppose $\tilde{T}_v \tilde{T}^{-1}_{w^{-1}}$ has a {\em minimal expression}, i.e., suppose we can write
$$
\tilde{T}_v \tilde{T}^{-1}_{w^{-1}} = \tilde{T}^{\epsilon_1}_{t_1} \cdots \tilde{T}^{\epsilon_k}_{t_k},
$$
where $t_1 \cdots t_k = vw$ is a reduced expression, and $\epsilon_i \in \{ 1,-1 \}$, for each $i$.  Letting $c_x$ denote the coefficient of $\tilde{T}_x$ in the expansion of such a minimal expression, one can easily prove by induction on $k$ that the $Q$-degree of $c_x$ is bounded above:
$$
{\rm deg}_Q  c_x \leq k - \ell(x).
$$ 

We get the following consequence, which we record for later use.

\begin{stz} \label{minimal_expression}
Suppose $\tilde{T}_v\tilde{T}^{-1}_{w^{-1}}$ has a minimal expression.  Then 
$$
{\rm deg}_Q R^{v}_{x,w}(Q) \leq \ell(vw) - \ell(x).
$$
\end{stz}

\subsection{The situation at $q=1$}

Recall from section \ref{hecke_algebra} that the elements 
 $C''_w = \Tr(\Fr_q, IC_w) = \Tr(\Fr_q, IC({\mathcal B}_{\bar{w}})[\ell(w)])$ form a basis of the Iwahori-Hecke
 algebra $\cH$. 
Consider the algebra homomorphism ${\mathcal H} \rightarrow {\Bbb Z}[\widetilde{W}]$ given by specializing $q^{1/2}=1$.  We will denote the image of $C''_w$ by $C''_w|$, etc.  The elements $T_w$ map to the usual basis of the group ring ${\Bbb Z}[\widetilde{W}]$, and the elements 
$C''_w$ map to another basis.  In fact
$$
C''_w| = \sum_{y \leq w} \varepsilon_wP_{y,w}(1) T_y|.
$$
We will identify $T_y|$ with $y$ in what follows.

Now fix a coweight $\lambda$, and write $\lambda = \lambda_1 - \lambda_2$, where both $\lambda_i$ are dominant.  We have 
$$
\varepsilon_\lambda q^{1/2}_\lambda {\Theta}_\lambda = \varepsilon_\lambda q^{1/2}_\lambda \widetilde{T}_{\lambda_1} \widetilde{T}^{-1}_{\lambda_2} = \sum_w a_w C''_w.
$$
Since the Wakimoto sheaf giving rise to this function is perverse and ${\mathcal B}$-equivariant, we see that $a_w \in \N[q^{1/2},q^{-1/2}]$ for every $w$, and thus $a_w \neq 0$ if and only if $a_w(1) \neq 0$.  We call the set $\lbrace w \,\, | \,\, a_w \neq 0 \rbrace$ the {\em IC-support} of $\varepsilon_\lambda q^{1/2}_\lambda {\Theta}_\lambda$.

Applying the homomorphism $x \mapsto x|$ to $\varepsilon_\lambda q^{1/2}_\lambda \widetilde{T}_{\lambda_1} 
\widetilde{T}^{-1}_{\lambda_2}$ gives 
$$
\varepsilon_\lambda t_{\lambda_1}t_{-\lambda_2} = \varepsilon_\lambda t_\lambda = 
\sum_{w \leq t_\lambda} Q_{w,t_\lambda}(1)C''_w|.
$$

We thus see that 
$$
a_{w}(1) = Q_{w,t_\lambda}(1),
$$
and that the IC-support of $\varepsilon_\lambda q^{1/2}_\lambda {\Theta}_\lambda$ is precisely the set $\lbrace w \leq t_\lambda \rbrace$.

Using now the geometric interpretation of inverse Kazhdan-Lusztig polynomials in the affine case (cf. Theorem 3.5), we get the following proposition.

\begin{stz}
The IC-support of the function $\varepsilon_\lambda q^{1/2}_\lambda {\Theta}_\lambda$ is precisely the set $\lbrace w \leq t_\lambda \rbrace$.  For each such $w$, the sum of the multiplicities $\sum_i m(\lambda,w,i)$ is equal to the intersection Euler characteristic of the space
$$
\cB_{t_\lambda}\cap \cB^{\overline{w}}.
$$
\end{stz}

\section{Applications to Shimura varieties}

\subsection{Relating nearby cycles on local models and Shimura varieties with Iwahori level structure}

In certain cases, the singularities occurring in the bad reduction of Shimura
varieties of PEL type with Iwahori level structure are equivalent to the
singularities of one of the schemes $M_\mu/\Z_p$ that we were considering. More
precisely, for minuscule $\mu$, $M_\mu$ is a local model in the sense of
Rapoport and Zink \cite{RZ2} of a certain Shimura variety.

Consider a Shimura variety $Sh = Sh_K(G,X)$ of PEL type, with Iwahori level
structure at $p$. (The same methods should yield analogous
results for arbitrary parahoric level structure.) Denote by $E$ its reflex
field. 
Since Shimura varieties of PEL type are moduli spaces of abelian varieties
with additional structure, there is a natural way to define a model over the
ring of integers of the completion of $E$ at a prime $\mathfrak p$ over $p$
(see \cite{Ko}, \cite{RZ2}). 

Furthermore, Rapoport and Zink have defined a so-called local model of the
Shimura variety: a scheme defined in terms of linear algebra, which has 'the
same' singularities as the model of the Shimura variety.

Now assume that $G_{\Q_p}$ is isomorphic to $GL_{n, \Q_p} \times \G_{m, \Q_p}$
or $GSp_{2g, \Q_p}$.

The Shimura datum gives rise to a minuscule cocharacter $\mu$ of $GL_{n,\Q_p}$
or $GSp_{2g,\Q_p}$, respectively. The functorial description of the local
model shows that it can be embedded into the deformation from $\Grass_{\Q_p}$
to $\Fl_{\F_p}$ \`a la Laumon (see \cite{HN1}), with generic fibre $\mathcal
Q_\mu$. Since the local model is flat \cite{G1}, \cite{G2}, it coincides with the scheme-theoretic
closure $M_\mu$ of $\mathcal Q_\mu$ in this deformation.

The relation between the model of the Shimura variety over $\Z_p$ and its
local model is given by a diagram 
$$\xymatrix{Sh & {\mathcal M} \ar[l]_\varphi \ar[r]^\psi & \Mloc
}$$
of $O_E$-schemes, where $\varphi$ is a torsor under a smooth affine group
scheme, and $\psi$ is smooth. The fibres of $\varphi$ are geometrically
connected (more precisely, this holds for the restriction of $\varphi$ to any 
geometric connected component of $\mathcal M$).  One can show that \'etale-locally around each
point of the special fibre, the model of the Shimura variety and the local
model are isomorphic. (See \cite{RZ2}.)

The stratification of the special fibre of $\Mloc$ induces stratifications of
the special fibres of $\mathcal M$ and $Sh$: the inverse images under $\psi$
of strata in $\Mloc_{\F_p}$ define a stratification of $\mathcal M$,
and there is a stratification of $Sh$ such that the inverse images of strata
under $\varphi$ are the strata in $\mathcal M$.  For details see \cite{GN}.  Following loc.cit., we call the resulting stratification on $Sh_{\F_p}$ the {\em Kottwitz-Rapoport stratification.}

The smooth base change theorem implies that the nearby cycles functor commutes
with pull-back under smooth morphisms, so we have (up to shifts by the dimensions)
$$ \varphi^\ast R\Psi_{Sh} \cong R\Psi_{\mathcal M} \cong \psi^\ast
R\Psi_{\Mloc}. $$
Furthermore, the intermediate extension commutes with pull-back under smooth
morphisms (see the proof of theorem III.11.3 in \cite{KW}), which implies that
$$ \varphi^\ast IC_{Sh, w} \cong IC_{{\mathcal M}, w} \cong \psi^\ast
IC_{\Mloc, w} \footnote{As in the introduction $IC_{Sh,w}$ and $IC_{\mathcal M, w}$ should be interpreted as direct sums of intersection complexes over geometric connected components of $Sh$ and $\mathcal M$.}.$$
(Again, this is up to a shift by the dimensions of the strata: these intersection complexes are normalized so that they are perverse.)

So we see that (since $\phi$ and $\psi$ have the same relative dimension) 
$$ \varphi^\ast R\Psi^{ss}_{Sh} 
  \cong \varphi^\ast \left(\bigoplus_w \bigoplus_i IC_{Sh, w}(-i)^{m(w,i)}\right),$$
and because $\varphi$ (restricted to any geometrically connected component of ${\mathcal M}$) is smooth with geometrically connected fibres,
this yields
$$ R\Psi^{ss}_{Sh} \cong \bigoplus_w \bigoplus_i IC_{Sh, w}(-i)^{m(w,i)}.$$

Hence the sheaf of nearby cycles on the Shimura variety decomposes in the same
way as the sheaf of nearby cycles on the local model.

\subsection{The weight spectral sequence}

By virtue of the Appendix, Theorem \ref{mixedness}, the nearby cycles $R\Psi := R\Psi_{Sh}$ is a {\em mixed} perverse sheaf on the special fiber of $Sh$.  Therefore, $R\Psi$ possesses a weight filtration (\cite{BBD}, $\S5$), that is, an increasing filtration $W$ such that each graded piece
$$
{\rm gr}^W_i := {\rm gr}^W_i(R\Psi)
$$
is a pure perverse sheaf of weight $i$.  Now assume that $Sh$ is proper over ${\mathcal O}_E$ (here $E$ denotes the completion of the reflex field at a prime above $p$).  Then we have a canonical isomorphism of Galois modules
$$
H^i(Sh_{\overline{\mathbb F}_p}, R\Psi) = H^i(Sh_{\overline{\mathbb Q}_p}, \Ql[\ell(\mu)]),
$$   
and the weight filtration on $R\Psi$ yields the {\em weight spectral sequence}
$$
_WE^{pq}_1 = H^{p+q}(Sh_{\overline{\mathbb F}_p}, gr^W_{-p}) \Longrightarrow 
H^{p+q}(Sh_{\overline{\mathbb Q}_p}, \Ql[\ell(\mu)]).
$$
This spectral sequence degenerates in $E_2$, and abuts to the weight filtration on the cohomology group on the right hand side.  Deligne's monodromy-weight conjecture asserts that the weight filtration on the right hand side is the same as the monodromy filtration defined with the help of the inertia action. 

One might hope that an explicit description of the $E_1$-terms might provide a modest first step toward the monodromy-weight conjecture for these Shimura varieties.  The methods of this paper allow us to completely describe the semi-simplification $({\rm gr}^W_i)^{ss}$ of each graded piece ${\rm gr}^W_i(R\Psi)$, in the category $P(Sh \times_s \eta, \Ql)$.  This can be used to describe the $E_1$-terms, if one forgets the Galois action.  The following summarizes our results in this direction.  (We will write $m(w,i)$ for $m(R\Psi,w,i)$ here.)

\begin{stz} \label{weight_ss}  In the category $P(Sh \times_s \eta, \Ql)$, we have the identity
$$
({\rm gr}^W_i)^{ss} = \bigoplus_{\gfrac{w \in {\rm Adm}(\mu)}{\ell(w) + 2j = i}} IC_{Sh,w}(-j)^{\oplus m(w,j)}.
$$
Here $IC_{Sh,w}$ is the intersection complex of the closure of the Kottwitz-Rapoport stratum $Sh_w$ in 
$Sh_{\overline{\F}_p}$.  Consequently, if $Sh$ is proper over $\mathcal O_E$, we deduce (forgetting the Galois actions) isomorphisms of $\Ql$-spaces
$$
_WE_1^{pq} = \bigoplus_{\gfrac{w \in {\rm Adm}(\mu)}{\ell(w) + 2j = -p}} 
IH^{p+q+\ell(w)}(\overline{Sh}_w, \Ql)^{\oplus m(w,j)}.
$$
\end{stz}

These results illustrate the fundamental importance of the Kottwitz-Rapoport stratification on $Sh_{\bar{\F}_p}$.  That stratification is indexed by the familiar set ${\rm Adm}(\mu)$, but the strata themselves are still rather mysterious.  In future work we will examine these strata more closely.

\section{Examples}

In this section we present a few examples for multiplicities
of nearby cycles sheaves $R\Psi = R\Psi(IC(\cQ_\mu)[\ell(t_\mu)])$. Most of the examples 
are for $G=GL_n$, but there is one example for each of $GSp_4$ and $GSp_6$, and also one for $G$ of
type $G_2$.

\subsection{The case of a divisor with normal crossings}

\begin{stz} In the Drinfeld case, where $GL_n$ and $\mu = (1, 0, \dots, 0)$,
we have $m(x,i) = 1$ for all $x\in\Adm(\mu)$, 
 $0\le i \le \ell(t_\mu)-\ell(x)$, and $m(x,i) = 0$ otherwise.
\end{stz}

{\em Proof.} This is a combinatorial exercise using the following facts.  In this case all Kazhdan-Lusztig polynomials are 1 (since the special fiber is a union of smooth divisors with normal crossings) 
and $\Tr(\Frob_q, R\Psi_x) = \varepsilon_\mu (1-q)^{\ell(t_\mu)-\ell(x)}$ (cf. \cite{H3}, Prop. 5.2).  \qed

\subsection{$GL_2$}

\begin{stz} In the case of $GL_2$ and an arbitrary coweight $\mu$,
all multiplicities $m(x,i) = 1$ for $x \in {\rm Adm}(\mu)$ and $0 \le i \le \ell(t_\mu) - \ell(x)$, and $m(x,i) = 0$ otherwise.
\end{stz}

{\em Proof.}  Argue as above, using the explicit formula for $\Tr(\Fr_q,R\Psi_x)$ given in \cite{H1}, Prop. 10.3. \qed 

\subsection{Other examples}

Further examples can be computed in the following way: the Kottwitz
conjecture (see section 2.7) describes the semi-simple trace of Frobenius on
$R\Psi$, and since the inertia action is unipotent, this is the same as the
usual trace. So the function of $R\Psi$ is known, and thus the
multiplicities can be computed by applying the base change from the
 $T_w$-basis to the $C''_w$-basis in the Hecke algebra. In order to carry
 this out, though, Bernstein functions and lots of Kazhdan-Lusztig
 polynomials have to be computed, so it is almost impossible to do this
 manually.
The following examples were computed with the help of a computer program. 

Because in all the interesting examples the number of alcoves
is quite large, we did not print out the multiplicities for
each single alcove, but we subsumed the information
for several alcoves into one line.

Besides the length $l$ and the number of alcoves that have certain
multiplicities, we also put the 'Bruhat configuration' in the tables,
i.e. the number of alcoves of length $l+1$, $l+2$, \dots
which are greater than the
given alcove with respect to the Bruhat order; cf. the remarks below.

The multiplicities listed below are always the numbers $m(w,i)$ where $i$ is understood to range from $0$ to 
$\ell(t_\mu) -\ell(w)$ (the multiplicities for $i$ not in this range always happen to vanish, cf. section 9.4 for some discussion of this point). 

\vskip.3cm



\subsubsection{$G=GL_4, \mu=(1,1,0,0)$}

Number of admissible alcoves: 33
\nopagebreak

\vskip.3cm

{\small\sf
\begin{tabular}{|l|l|c|l|}
\hline
Length & \#Alcoves & Multiplicities & Bruhat configuration \\
\hline \hline

l=0 & 1 & 1,1,2,1,1 & 4, 10, 12, 6 \\ \hline
l=1 & 4 & 1,2,2,1 & 5, 8, 5\\ \hline
l=2 & 8 & 1,1,1 & 3, 3\\
 & 2 & 1,2,1 & 4, 4 \\ \hline
l=3 & 12 & 1,1 & 2 \\ \hline
l=4 & 6 & 1 & - \\ \hline
\end{tabular}
}

This is the simplest non-trivial case. We see already that the
multiplicities for $w$ are not determined by the length of $w$.

\subsubsection{$G=GL_5, \mu=(1,1,0,0,0)$}

Number of admissible alcoves: 131
\nopagebreak

\vskip.3cm

{\small\sf
\begin{tabular}{|l|l|c|l|}
\hline
Length & \#Alcoves & Multiplicities & Bruhat configuration \\
\hline \hline

l=0 & 1 & 1, 1, 2, 2, 2, 1, 1 & 5, 15, 30, 40, 30, 10 \\ \hline
l=1 & 5 & 1, 2, 3, 3, 2, 1 &  6, 17, 28, 24, 9 \\ \hline
l=2 & 5 & 1, 1, 2, 1, 1 & 4, 10, 12, 6 \\
    & 5 & 1, 2, 2, 2, 1 & 6, 14, 15, 7 \\
    & 5 & 1, 3, 4, 3, 1 & 7, 17, 18, 8 \\ \hline
l=3 & 5 & 1, 1, 1, 1 & 4, 6, 4 \\
    & 20 & 1, 2, 2, 1 & 5, 8, 5 \\
    & 5 & 1, 2, 2, 1 & 5, 9, 6\\ \hline
l=4 & 30 & 1, 1, 1 & 3, 3\\
    & 10 & 1, 2, 1 & 4, 4\\ \hline
l=5 & 30 & 1, 1 & 2\\ \hline
l=6 & 10 & 1 & - \\ \hline
\end{tabular}
}

\subsubsection{$G=GL_6, \mu=(1,1,0,0,0,0)$}

Number of  admissible alcoves: 473
\nopagebreak

\vskip.3cm

{\small\sf
\begin{tabular}{|l|l|c|l|}
\hline
Length & \#Alcoves & Multiplicities & Bruhat configuration \\
\hline \hline

l=0 & 1 & 1, 1, 2, 2, 3, 2, 2, 1, 1 & 6, 21, 50, 90, 120, 110, 60, 15 \\ \hline
l=1 & 6 & 1, 2, 3, 4, 4, 3, 2, 1 & 7, 24, 55, 86, 88, 52, 14 \\ \hline
l=2 & 6 & 1, 1, 2, 2, 2, 1, 1 & 5, 15, 30, 40, 30, 10 \\
    & 6 & 1, 2, 3, 3, 3, 2, 1 & 7, 24, 48, 58, 39, 12 \\
    & 9 & 1, 3, 5, 6, 5, 3, 1 & 8, 28, 56, 67, 44, 13 \\ \hline
l=3 & 30 & 1, 2, 3, 3, 2, 1 & 6, 17, 28, 24, 9 \\
    & 6 & 1, 2, 2, 2, 2, 1 & 7, 20, 30, 24, 9 \\
    & 12 & 1, 3, 4, 4, 3, 1 & 8, 25, 39, 31, 11 \\
    & 2 & 1, 4, 7, 7, 4, 1 & 9, 30, 47, 36, 12 \\ \hline
l=4 & 15 & 1, 1, 2, 1, 1 & 4, 10, 12, 6 \\
    & 6 & 1, 1, 1, 1, 1 & 5, 10, 10, 5 \\
    & 30 & 1, 2, 2, 2, 1 & 6, 14, 15, 7 \\
    & 6 & 1, 2, 2, 2, 1 & 6, 14, 16, 8 \\
    & 3 & 1, 2, 3, 2, 1 & 6, 15, 18, 9 \\
    & 30 & 1, 3, 4, 3, 1 & 7, 17, 18, 8 \\ \hline
l=5 & 30 & 1, 1, 1, 1 & 4, 6, 4 \\
    & 60 & 1, 2, 2, 1 & 5, 8, 5 \\
    & 30 & 1, 2, 2, 1 & 5, 9, 6 \\ \hline
l=6 & 80 & 1,1,1 &  3, 3 \\
    & 30 & 1,2,1 & 4, 4 \\ \hline
l=7 & 60 & 1,1 & 2 \\ \hline
l=8 & 15 & 1 & - \\ \hline
\end{tabular}
}

\subsubsection{$G=GL_3, \mu=(2,2,0)$}

Number of admissible alcoves: 19
\nopagebreak

\vskip.3cm

{\small\sf
\begin{tabular}{|l|l|c|l|}
\hline
Length & \#Alcoves & Multiplicities & Bruhat configuration \\
\hline \hline

l=0 & 1 & 1, 1, 2, 1, 1 & 3, 6, 6, 3 \\ \hline
l=1 & 3 & 1,2,2,1 &  4, 5, 3 \\ \hline
l=2 & 3 & 1,1,1 &  2, 2 \\
    & 3 & 1,1,1 & 3, 3 \\ \hline
l=3 & 6 & 1,1 & 2 \\ \hline
l=4 & 3 & 1 & - \\ \hline
\end{tabular}
}

\subsubsection{$G=GL_3, \mu=(3,1,0)$}

Number of admissible alcoves: 49
\nopagebreak

\vskip.3cm

{\small\sf
\begin{tabular}{|l|l|c|l|}
\hline
Length & \#Alcoves & Multiplicities & Bruhat configuration \\
\hline \hline

l=0 & 1 & 1, 2, 3, 3, 3, 2, 1 & 3, 6, 9, 12, 12, 6 \\ \hline
l=1 & 3 & 1, 3, 5, 5, 3, 1 & 4, 8, 12, 12, 6 \\ \hline
l=2 & 3 & 1, 2, 3, 2, 1 & 4, 9, 11, 6 \\ 
    & 3 & 1, 3, 4, 3, 1 & 4, 9, 11, 6 \\ \hline
l=3 & 3 & 1,2,2,1 & 4, 6, 4 \\
    & 3 & 1,2,2,1 & 4, 8, 6 \\
    & 3 & 1,4,4,1 & 6, 10, 6 \\ \hline
l=4 & 3 & 1,1,1 & 2, 2 \\
    & 9 & 1,2,1 & 4, 4 \\ \hline
l=5 & 12 & 1,1 & 2 \\ \hline
l=6 & 6 & 1 & - \\ \hline
\end{tabular}
}

This example shows that the Bruhat configuration {\em numbers} given in the
table do not determine the multiplicities (look at the elements of length
2). On the other hand, in this example it is still quite easy to check
that the Bruhat {\em graphs} for the two types of elements are different.

\subsubsection{$G=GL_4, \mu=(2,0,0,0)$}

Number of admissible alcoves: 65
\nopagebreak

\vskip.3cm

{\small\sf
\begin{tabular}{|l|l|c|l|}
\hline
Length & \#Alcoves & Multiplicities & Bruhat configuration \\
\hline \hline

l=0 & 1 & 1, 1, 2, 2, 2, 1, 1 & 4, 10, 16, 18, 12, 4 \\ \hline
l=1 & 4 & 1, 2, 3, 3, 2, 1 & 5, 11, 15, 11, 4\\ \hline
l=2 & 4 & 1, 1, 2, 1, 1 & 3, 6, 6, 3 \\ 
    & 4 & 1, 2, 2, 2, 1 & 5, 10, 9, 4\\
    & 2 & 1, 3, 4, 3, 1 & 6, 12, 10, 4 \\ \hline
l=3 & 12 & 1,2,2,1 & 4, 5, 3 \\
    & 4 & 1,1,1,1 & 4, 6, 4 \\ \hline
l=4 & 6 & 1,1,1 & 2, 2 \\
    & 12 & 1,1,1 & 3, 3 \\ \hline
l=5 & 12 & 1,1 & 2\\ \hline
l=6 & 4 & 1 & - \\ \hline
\end{tabular}
}

\subsubsection{$G=GL_4, \mu=(2,1,0,0)$}

Number of admissible alcoves: 143
\nopagebreak

\vskip.3cm

{\small\sf
\begin{tabular}{|l|l|c|l|}
\hline
Length & \#Alcoves & Multiplicities & Bruhat configuration \\
\hline \hline

l=0 & 1 & 1, 2, 3, 4, 4, 3, 2, 1 & 4, 10, 20, 30, 36, 30, 12 \\ \hline
l=1 & 4 & 1, 3, 5, 6, 5, 3, 1 & 5, 14, 25, 33, 29, 12 \\ \hline
l=2 & 4 & 1, 3, 4, 4, 3, 1 & 5, 12, 20, 21, 10 \\
    & 4 & 1, 2, 3, 3, 2, 1 & 5, 14, 25, 26, 12 \\
    & 2 & 1, 6, 11, 11, 6, 1 & 8, 20, 30, 28, 12\\ \hline
l=3 & 4 & 1,2,2,2,1 & 3, 6, 9, 6\\
    & 4 & 1,2,3,2,1 & 5, 13, 17, 9 \\
    & 12 & 1,3,4,3,1 & 6, 15, 19, 10\\ \hline
l=4 & 12 & 1,2,2,1 & 4, 8, 6 \\
    & 8 & 1,2,2,1 & 5, 9, 6\\
    & 8 & 1,3,3,1 & 6, 11, 7\\
    & 2 & 1,3,3,1 & 6, 12, 8 \\ \hline
l=5 & 8 & 1,1,1 & 3, 3 \\
    & 28 & 1,2,1 & 4, 4 \\ \hline
l=6 & 30 & 1,1 & 2 \\ \hline
l=7 & 12 & 1 & - \\ \hline
\end{tabular}
}

\subsubsection{$G=GSp_4$, $\mu = (1,1,0,0)$}

Number of admissible alcoves: 13
\nopagebreak

{\small\sf
\begin{tabular}{|l|l|c|l|}
\hline
Length & \#Alcoves & Multiplicities & Bruhat configuration \\
\hline \hline
l=0 & 1 & 1, 1, 1, 1 & 3, 5, 4 \\ \hline
l=1 & 2 & 1, 1, 1 & 3, 3 \\ 
    & 1 & 1, 2, 1 & 4, 4 \\ \hline
l=2 & 5 & 1,1 & 2 \\ \hline
l=3 & 4 & 1 & - \\ \hline
\end{tabular}
}

\subsubsection{$G=GSp_6$, $\mu = (1,1,1,0,0,0)$}

Number of admissible alcoves: 79
\nopagebreak

\vskip.3cm

{\small\sf
\begin{tabular}{|l|l|c|l|}
\hline
Length & \#Alcoves & Multiplicities & Bruhat configuration \\
\hline \hline
l=0 & 1 & 1,1,1,2,1,1,1 & 4, 9, 17, 22, 18, 8 \\ \hline
l=1 & 2 & 1,1,2,2,1,1 & 4, 10, 16, 15, 7 \\
    & 2 & 1,2,3,3,2,1 & 5, 13, 19, 17, 8 \\ \hline
l=2 & 1 & 1,1,2,1,1 & 4, 10, 12, 6 \\
    & 4 & 1,2,2,2,1 & 5, 10, 11, 6 \\
    & 2 & 1,2,2,2,1 & 5, 11, 14, 8 \\
    & 2 & 1,3,4,3,1 & 6, 13, 14, 7 \\ \hline
l=3 & 6 & 1,1,1,1 & 3, 5, 4 \\
    & 2 & 1,1,1,1 & 4, 6, 4 \\
    & 4 & 1,2,2,1 & 5, 8, 5 \\
    & 4 & 1,2,2,1 & 5, 9, 6 \\
    & 1 & 1,3,3,1 & 6, 12, 8 \\ \hline
l=4 & 14 & 1,1,1 & 3, 3 \\
    & 8 & 1,2,1 & 4, 4 \\ \hline
l=5 & 18 & 1,1 & 2 \\ \hline
l=6 & 8 & 1 & - \\ \hline
\end{tabular}
}

\subsubsection{$G$ of type $G_2$, $\mu = (2,1,0$)}

Number of admissible alcoves: 41
\nopagebreak

\vskip.3cm

{\small\sf
\begin{tabular}{|l|l|c|l|}
\hline
Length & \#Alcoves & Multiplicities & Bruhat configuration \\
\hline \hline

l=0 & 1 & 1, 1, 1, 1, 1, 1, 1 & 3, 5, 7, 9, 10, 6 \\ \hline
l=1 & 1 & 1, 1, 1, 1, 1, 1 & 3, 5, 7, 9, 6 \\
    & 1 & 1, 1, 1, 1, 1, 1 & 3, 6, 9, 10, 6 \\
    & 1 & 1, 2, 2, 2, 2, 1 & 4, 7, 9, 10, 6 \\ \hline
l=2 & 2 & 1,1,1,1,1 & 3, 5, 7, 5 \\
    & 1 & 1,2,2,2,1 & 4, 7, 9, 6 \\
    & 2 & 1,2,2,2,1 & 4, 8, 10, 6 \\ \hline
l=3 & 3 & 1,1,1,1 & 3, 5, 4 \\
    & 2 & 1,2,2,1 & 4, 7, 5 \\
    & 1 & 1,2,2,1 & 4, 8, 6 \\
    & 1 & 1,2,2,1 & 5, 9, 6 \\ \hline
l=4 & 4 & 1,1,1 & 3, 3 \\
    & 5 & 1,2,1 & 4, 4 \\ \hline
l=5 & 10 & 1,1 & 2 \\ \hline
l=6 & 6 & 1 & - \\ \hline
\end{tabular}
}

\subsection{Numerical Observations and Conjectures}

Looking at the tables, one can make several empirical observations:

\begin{enumerate}
\item[(A)]  We have $m(w,i) \in [0, \ell(\mu) - \ell(w)]$, or equivalently, the function ${\rm Tr}({\rm Fr}_q, R\Psi_w)$ is a polynomial in $q$ of degree $\leq \ell(\mu) - \ell(w)$;
\item[(B)]  The sequence $m(w,0), \dots, m(w, \ell(\mu) - \ell(w))$ is palindromic;  it increases to the middle, then decreases again.
\item[(C)]  We have $m(w,0) = m(w, \ell(\mu) - \ell(w)) = 1$.
\end{enumerate}

From some more detailed information not contained in the tables, we also made the following striking observation:

\noindent $\bullet$   The multiplicity polynomial $m(w) = \sum_i m(w,i) q^i$ is determined by the isomorphism type of the Bruhat graph of the set $\{ x \in {\rm Adm}(\mu) ~ | ~ x \geq w \}$ (although not by the Bruhat configuration numbers alone).  It is remarkable that the multiplicity polynomials seem to be determined by data that is independent of the underlying root 
system.

\medskip

Using the explicit formula for ${\rm Tr}({\rm Fr}_q, R\Psi_w)$ in the minuscule case, proved in \cite{H3},\cite{HP}, we can prove this last observation in that case.  Based on our data we conjecture it holds in complete generality.  We offer the following explanations for (A)-(C).

\medskip

\noindent {\em Proof of (A) for $GL_n$, or for $\mu$ minuscule.} 

From the definition, we get a recursive procedure to compute the multiplicities:

$$\varepsilon_w \sum_i m(w,i)q^i = \Tr(\Fr_q, R\Psi_w) - \sum_{\gfrac{x\in \Adm(\mu)}{x >
    w}} \varepsilon_x m(x) P_{w,x}.$$

By downward induction on $\ell(w)$, we deduce that ${\rm Tr}({\rm Fr}_q, R\Psi_w)$ is a polynomial of degree $\leq \ell(\mu) - \ell(w)$ for all $w \in {\rm Adm}(\mu)$ if and only if $m(w)$ is a polynomial of degree 
$\leq \ell(\mu) - \ell(w)$, for all such $w$.

Now suppose $\mu$ is either a minuscule coweight, or an arbitrary coweight for the group $GL_n$.  In \cite{H3} and \cite{HP} it is proved that for any coweight $\nu \in \Omega(\mu)$, the function $\Theta_\nu$ has a minimal expression.  Therefore by Proposition \ref{minimal_expression}, we see that
$$
\Theta_\nu(w) = q_w^{-1/2}R(Q),
$$
where $R(Q) \in {\mathbb N}[Q]$, and ${\rm deg}_QR(Q) \leq \ell(t_\nu) - \ell(w)$.  From this it is easy to see (via the Kottwitz conjecture), that ${\rm Tr}({\rm Fr}_q, R\Psi_w)$ is a polynomial in $q$ of degree $\leq \ell(\mu) - \ell(w)$.

We remark that the conclusion also holds when every dominant coweight 
$\lambda \leq \mu$ is a sum of dominant minuscule coweights.

\medskip

\noindent {\em Proof of (B).}

Recall that $DR\Psi = R\Psi(\ell(\mu))$.  Using Lemma \ref{proto_palindromic}, we deduce immediately that the multiplicity function $m(R\Psi,w) := \sum_i m(R\Psi,w,i) q^i$ is {\em palindromic}, i.e., 
$\overline{m(R\Psi,w)} = q^{-1}_{\mu} q_w m(R\Psi,w)$.

The multiplicities being palindromic is also a consequence of the
conjecture that the monodromy and the weight filtration on nearby cycles agree. In 
the function field case this is known to be true (cf. \cite{BB} $\S 5$); since the multiplicities
in the unequal characteristic case and in the function field case coincide,
we get the same consequences in the unequal characteristic case.
In particular this proves also that the multiplicities increase to the 
middle, and then decrease again.

\medskip

\noindent {\em Proof of (C) in the minuscule case.}

It is easy to see that, for minuscule $\mu$,
 $ m(w, l(\mu)-l(w))=1$ for all $w$ because then $\Tr(\Fr_q, R\Psi_w)$ 
is essentially an $R$-polynomial (\cite{H3}, \cite{HP}), and thus the coefficient 
of $q^{l(\mu)-l(w)}$ in $\Tr(\Fr, R\Psi_w)$ is  $\varepsilon_w$. 

This also implies $m(w,0)=1$  for all $w$, since the
multiplicities are palindromic.

If $\mu$ is minuscule, the identity $m(x,0)=1$ for all $x \in {\rm Adm}(\mu)$ is equivalent to the following combinatorial identity for such $x$:
$$
\sum_{\gfrac{w\in \Adm(\mu)}{w \ge x}} \varepsilon_w = \varepsilon_\mu.
$$ 

This has certain geometric consequences.  For example, it shows that any codimension 1 Iwahori-orbit in a local model is contained in at least two irreducible components.  It follows from this and properties of the Bruhat order that the smooth locus of a local model is precisely the union of the extreme strata.  Consequently, we get the analogous results for certain Shimura varieties.  (See \cite{GN} for related results, proved by a different method.)

\begin{stz} Let $Sh$ be a Shimura variety as in section 8.  Then the smooth locus of the special fiber of 
$Sh$ is precisely the union of the Kottwitz-Rapoport strata indexed by the translation elements in the set ${\rm Adm}(\mu)$.  Here $\mu$ is the minuscule coweight in the Shimura datum for $Sh$.
\end{stz}

\pagebreak

\section{Appendix}

\begin{center}
\Large
Part I: Nearby cycles preserve mixedness
\end{center}

Due to fundamental results of Deligne \cite{Weil2}, it is known that nearby cycles preserve mixedness in the function-field setting (i.e., when working over a trait which is the Henselization at a point of a smooth 
curve over a finite field).  The aim of this appendix is to prove an analogous result in the setting of unequal characteristic traits.  This result applies to nearby cycles on local models and integral models for Shimura varieties over $p$-adic number rings, ensuring that they possess weight filtrations, as postulated in the introduction.

\medskip

Let $(S,s,\eta)$ denote a Henselian trait (i.e, the spectrum of a complete 
discrete valuation ring).  We assume that the residue field $k(s)$ is finite.  
Let $k(\bar{\eta})$ be a separable closure of $k(\eta)$, 
and let $\bar{S}$ denote the normalization of $S$ in $\bar{\eta}$, with closed point $\bar{s}$.  
Let $X \rightarrow S$ be a separated, finite-type scheme over 
$S$.  

\medskip

We will use the notion of $S$-variety (cf. \cite{deJ2}, 2.15).  A separated finite-type $S$-scheme $X$ will be called an $S$-variety if $X \rightarrow S$ is 
flat and $X$ is integral.

\medskip

For $A \in D^b_c(X_\eta,\bar{\mathbb Q}_\ell)$, define 
$R\Psi^X(A) = \bar{i}^*R\bar{j}_*
A_{\bar{\eta}}$, an object in $D^b_c(X \times_s \eta, \bar{\mathbb Q}_\ell)$.

\medskip

Let $\Lambda_\eta$ denote the constant $\ell$-adic sheaf $\, \bar{\mathbb Q}_\ell$ on 
$X_\eta$ (or some $\lambda$-adic constant sheaf $E_\lambda$, for $\lambda | \ell$) 
We will prove the following theorem. 

\medskip

\begin{thm}\label{mixedness}
Suppose that $X_\eta$ is smooth.  
Then the complex $R\Psi^X(\Lambda_\eta)$ is {\em mixed}.
\end{thm}
  
\medskip

The proof of the theorem consists of several steps which eventually reduce it to the calculations of Rapoport-Zink (\cite{RZ1}, \cite{I}) for the case of a proper strictly semi-stable $S$-variety (for the precise definition of {\em strictly semi-stable}, see \cite{deJ2}, 2.16).  We will use de Jong's alterations \cite{deJ2} to reduce the general case to that case.  Roughly, if $X/S$ is a proper $S$-variety, one takes a strictly semi-stable alteration $\pi: X' \rightarrow X$, and then applies the push-forward functor $\pi_*$ to appropriate nearby cycles on the alteration.  Thanks to Rapoport-Zink and properties of the push-forward functor, the resulting complex is mixed.  
One proves using the trace map that it contains the original nearby cycles as a subobject, thus essentially completing the proof of the theorem in that case.  (An extra difficulty arises because in the end we have to descend from a finite extension $k(s')$ back down 
to $k(s)$; see Corollary \ref{descend}.)  It is convenient to work with perverse sheaves throughout this process.  In particular, we will work with the perverse push-forward $^p\pi_*$.  

\medskip

\subsection{Alterations}

The following theorem is a direct consequence of \cite{deJ2}, Theorem 6.5.

\begin{thm} [de Jong] \label{alterations}
Let $X/S$ be an $S$-variety.  Then there exists a trait $S' = (S',s',\eta')$ finite over $S$, an $S'$-variety $X'$ for which there is an alteration $X' \rightarrow X$ of $S$-varieties, and an open immersion $j: X' \hookrightarrow \overline{X'}$ of $S'$-varieties such that $\overline{X'}$ is a proper strictly semi-stable $S'$-variety:

$$\xymatrix{
\overline{X'} \ar[d] & X' \ar[l]_j \ar[d] \ar[r] & X \ar[d] \\ 
S' & S' \ar[l]_{\rm id} \ar[r] & S.}
$$
\end{thm}

\bigskip

\subsection{Proof for $X'/S'$}

Let $X', S'$ etc. be as in Theorem \ref{alterations}.  
The special fiber $\overline{X'}_{s'}$ is globally a union of reduced divisors with normal crossings and the calculations of Rapoport-Zink are valid 
(cf. \cite{I}, Th. 3.2).  The sheaves   
$R^q\Psi^{\overline{X'}}(\Lambda_{\eta'})$ are explicitly computed in loc. cit, and are easily seen to be mixed.  Thus Theorem \ref{mixedness} holds in this case. 

\medskip

This implies the theorem for $X'/S'$.  By Theorem \ref{alterations} the morphism $j$ embeds $X'$  as an open subscheme in $\overline{X'}$.  Then the formula 
$$
j^* R\Psi^{\overline{X'}}(\Lambda_{\eta'}) = R\Psi^{X'}(\Lambda_{\eta'})
$$
shows that the right hand side is mixed.

\medskip

\subsection{Proof for $X_\eta$ smooth and geometrically integral}

Assume that $X/S$ is finite-type and separated, with $X_\eta$ smooth and geometrically integral.  The nearby cycles are supported on the scheme-theoretic closure of $X_\eta$ in $X_s$.  We may therefore assume $X/S$ is flat by replacing $X$ with the scheme-theoretic closure of $X_\eta$ in $X$ (which is flat with the same generic fiber since $X_\eta$ 
is reduced).  Thus, we may assume $X/S$ is an $S$-variety whose generic fiber is smooth and geometrically integral.

\medskip

Let $S'$ and $X'$ be as in Theorem \ref{alterations}. 

\medskip

Let us denote by $\pi: X' \rightarrow X_{S'}$ the resulting alteration of $S'$-varieties.
Let us also denote by $\pi: X'_{\eta'} \rightarrow X_{\eta'}$ the morphism on the generic fibers.  This is an alteration, hence by definition is generically finite, say of degree $n$.  Since $X_{\eta'}$ is integral, $\pi$ is also generically flat.  More precisely, there exists a non-empty open subset $U \subset X_{\eta'}$ such that $\pi: \pi^{-1}(U) \rightarrow U$ is flat and finite.  
\medskip

Let us denote $^p\pi^*:= \, ^pH^0\pi^*$ (resp. $^p\pi_* := \, ^pH^0\pi_*$), the perverse versions of the pull-back and push-forward functors via $\pi$.  Since $\pi$ is proper, we have $\pi_* = \pi_!$ and so $^p\pi_* = \, ^p\pi_!$.  

\medskip

Let $\Lambda^{X}_{\eta'}$ denote the constant sheaf on $X_{\eta'}$, and let $\Lambda^{X'}_{\eta'}$ denote the constant sheaf on $X'_{\eta'}$.  For $d = {\rm dim}(X_\eta)$, let $A^{X} = 
\Lambda^{X}_{\eta'}[d](\frac{d}{2})$, and $A^{X'} = \Lambda^{X'}_{\eta'}[d](\frac{d}{2})$; the smoothness of $X'_{\eta'}$ (resp. $X_{\eta'}$) ensures that $A^{X'}$ (resp. $A^{X}$) is perverse of weight zero.  
We will need the following lemma.

\medskip

\begin{lem} \label{subobject} The perverse sheaf $A^X$ is a subobject of $\,\, 
^p\pi_*(A^{X'})$, in the category $P(X_{\eta'}, \Ql)$.  
\end{lem}

\noindent {\em Proof.}
The adjunction morphism $A^X \rightarrow \pi_* \pi^* A^X$ and the perversity of $A^X$ yield a morphism of perverse sheaves $A^X \rightarrow \, ^p\pi_* A^{X'}$.  Since $A^X$ is simple, it is sufficient to prove this morphism is not the zero morphism.  This can be detected by showing that for some dense open $U \subset X_{\eta'}$, the restriction of the morphism to $U$ is nonzero.  Now choose $U$ over which $\pi$ is flat and finite, and use the trace map 
$$
{\rm Tr} : (\pi_!\pi^*(\Lambda^{X}_{\eta'}))|_U \rightarrow (\Lambda^{X}_{\eta'})|_U,
$$
which composed with the canonical adjunction map $\Lambda^X_{\eta'} \rightarrow 
\pi_*\pi^*(\Lambda^X_{\eta'}) = \pi_! \pi^*(\Lambda^X_{\eta'})$ is just multiplication by the degree $n$.  For the definition and basic properties of this trace map, see [SGA 4, Exp. XVII, p. 553-554].
\qed

\bigskip

Now let $r: X_{S'} \rightarrow X$ denote the projection morphism.  
Since $r$ arises from a finite ``change of trait'' morphism, the invariance of nearby cycles under change of trait ([SGA 4 1/2, Thm. finitude 3.7]) implies that we have a canonical isomorphism $r^*R\Psi^X(\Lambda_\eta[d](\frac{d}{2})) = R\Psi^{X_{S'}}(A^{X})$ in the category 
$D^b_c(X \times_s \eta', \bar{\mathbb Q}_\ell)$.   Here the meaning of this last category and of the left hand side of the equality is the following.  The nearby cycles $R\Psi^X(\Lambda_\eta[d](\frac{d}{2}))$ live on the geometric special fiber $X_{\bar{s}} = X_{\bar{s'}}$, and naturally carry a compatible action of ${\rm Gal}(\bar{\eta}/\eta)$.  The sheaf 
$r^*R\Psi^X(\Lambda_\eta[d](\frac{d}{2}))$ is just that same sheaf together with the compatible action of ${\rm Gal}(\bar{\eta'}/\eta')$ via the canonical restriction map ${\rm Gal}(\bar{\eta'}/\eta') \rightarrow {\rm Gal}(\bar{\eta}/\eta)$.  By 
Lemma \ref{subobject}, this is a subobject (in the category $P(X \times_{s'} \eta', \Ql)$) of 
$$
R\Psi^{X_{S'}}(\, ^p\pi_*(A^{X'})) = \, ^p\pi_* R\Psi^{X'}(A^{X'}).
$$
(The last equality follows from the compatibility of nearby cycles and proper push-forwards together with the fact that $R\Psi$ commutes with the perverse truncation functors $^p\tau_{\leq 0}$ and $^p\tau_{\geq 0}$, hence with the functor $^pH^0$.  The compatibility of $R\Psi$ and the perverse truncation functors follows easily from the fact that $R\Psi$ preserves $^pD_c^{\leq 0}$ and commutes with Verdier duality.)  

\medskip
We have used here the basic fact that the functor $R\Psi$, being $t$-exact, induces an exact functor on categories of perverse sheaves.

\medskip

We can now complete the argument.
By the case of $X'/S'$ treated above, the complex $R\Psi^{X'}(A^{X'})$ is mixed.  Thus by  \cite{BBD} $\S 5$, $^p\pi_*R\Psi^{X'}(A^{X'})$ is mixed.  Its subobject $r^*R\Psi^X(\Lambda_\eta[d](\frac{d}{2}))$ is mixed too.  By the lemma and corollary in the following section on descent, we deduce finally that $R\Psi^X(\Lambda_{\eta})$ is also 
mixed.  This proves Theorem \ref{mixedness} in the case of finite-type separated $S$-schemes $X$ where $X_\eta$ is smooth and geometrically integral.
\medskip

\subsection{Proof for $X_\eta$ smooth}

We consider finally the general case, where $X_\eta$ is only assumed smooth.  Since $X_{\bar{\eta}}$ is smooth, its connected components are integral.  If $\eta'/\eta$ is a finite extension over which all the connected components are defined, then it follows that each connected component of $X_{\eta'}$ is smooth and geometrically integral.  Let $S'$ denote the normalization of $S$ in $\eta'$; then $S'$ is a trait, finite over $S$.  As above let $r: X_{S'} \rightarrow X$ denote the projection.  For each connected component 
$X^{(\alpha)}_{\eta'}$, $\alpha = 1, 2, \dots, m$, let $X^{(\alpha)}$ denote its scheme-theoretic closure in $X_{S'}$.  We have closed immersions of $S'$-schemes $i_{\alpha}: X^{(\alpha)} \hookrightarrow 
X_{S'}$.  Each $X^{(\alpha)}$ is a separated finite-type $S'$-scheme whose generic fiber is smooth and geometrically integral; hence the results of the previous section obtain for it.

It is clear that
$$
R\Psi^{X_{S'}}(\Lambda_{\eta'}) = \bigoplus_{\alpha} i_{\alpha, *} R\Psi^{X^{(\alpha)}}(\Lambda_{\eta'});
$$
by the previous section the right hand side is mixed, hence so is the left hand side.  Further, as before, $r^* R\Psi^{X}(\Lambda_\eta) = R\Psi^{X_{S'}}(\Lambda_{\eta'})$ is just the complex $R\Psi^X(\Lambda_\eta)$ endowed with its compatible action of ${\rm Gal}(\bar{\eta'}/\eta')$ via ${\rm Gal}(\bar{\eta'}/\eta') \rightarrow 
{\rm Gal}(\bar{\eta}/\eta)$.  Now Corollary \ref{descend} finally proves that $R\Psi^{X}(\Lambda_\eta)$ is mixed, as desired.

\medskip

\subsection{Descent step}

 The following results are well-known (cf. \cite{KW}, p.14); we include them for the reader's convenience.

\medskip

\begin{lem}\label{finite_surjective}
Let $r: S' \rightarrow S$ denote a finite surjective morphism of traits with finite residue fields.  Let $X/S$ be finite-type, and by abuse denote also by $r$ the morphism $X_{s'} \rightarrow X_s$ on special fibers deduced from $r \times_S {\rm id}_X: X_{S'} \rightarrow X$. 

If ${\mathcal G}$ is a smooth \'{e}tale Weil sheaf on $X_{s}$ such that $r^*{\mathcal G}$ is a mixed sheaf on $X'_{s'}$ with weights $\leq w$, then ${\mathcal G}$ is also mixed with weights $\leq w$.
\end{lem}

We suppose $k(s)$ has $q$ elements.
As noted above in section 10.3, because $X_{\bar{s}} = X_{\bar{s'}}$, we can regard ${\mathcal G}$ as a sheaf on $X_{\bar{s}}$ endowed with a compatible action of a geometric Frobenius ${\rm Frob}_q \in 
{\rm Gal}(\bar{\eta}/\eta)$, and we can regard $r^*{\mathcal G}$ as this same sheaf endowed with a compatible action of a geometric Frobenius in ${\rm Gal}(\bar{\eta'}/\eta')$ which maps to ${\rm Frob}^n_q = {\rm Frob}_{q^n}$ under the restriction map ${\rm Gal}(\bar{\eta'}/\eta') \rightarrow 
{\rm Gal}(\bar{\eta}/\eta)$, where $n = [k(s'):k(s)]$.  

\noindent {\em Proof.}
 Let $X_0 = X_s = s \times _S X$.  Note that the ``change of trait'' morphism $r$ above is the morphism 
$$
s' \times_{S'} X' \rightarrow s \times_S X
$$
which arises via the canonical identifications $s' \times_{S'} (S' \times_S X) = s' \times_S X = 
s' \times_s X_s$ from the projection
$$
r: X_0 \times_{\kappa} \kappa' \rightarrow X_0,
$$
where $X_0 = X_s$, $\kappa = {\mathbb F}_q$, and $\kappa' = {\mathbb F}_{q^n}$.

The sheaf $r^*{\mathcal G}$ is by hypothesis a smooth mixed Weil sheaf on $X_0 \times_\kappa \kappa'$.  We want to prove that ${\mathcal G}$ is a mixed Weil sheaf on $X_0$. 
  
Since $r^*$ is exact, we can assume without loss of generality that ${\mathcal G}$ is simple.  We will prove that in this case ${\mathcal G}$ is pure.  

Since $r^*{\mathcal G}$ is mixed, there is a non-zero smooth Weil sheaf ${\mathcal F}$ on $X_0 \times_\kappa \kappa'$ which is a simple (thus pure) Weil subsheaf of $r^*{\mathcal G}$.

Let now $F$ denote the Frobenius automorphism ${\rm Fr}_q$ for $X_0$, so that $F^n$ is the Frobenius automorphism for $X_0 \times_{\kappa} \kappa'$.  We are given an isomorphism $(F^{n})^*{\mathcal F} \tilde{\lto} {\mathcal F}$.  Using this and the fact that the category of Weil sheaves is closed under the formation of kernels and cokernels, one can show that 
$$
\sum_{i=0}^{n-1} (F^i)^* {\mathcal F}
$$
is a nonzero Weil subsheaf of ${\mathcal G}$.  Since the latter is simple, we have in fact
$$
\bigoplus_{i \in I} (F^i)^* {\mathcal F} = {\mathcal G},
$$
where $I \subset \{0,1,\dots, n-1 \}$ is a nonempty subset.  (The sum can be made a direct sum by omitting certain $i$'s.)

Now fix any isomorphism $\tau: \bar{\mathbb Q}_\ell \tilde{\lto} {\mathbb C}$.  Note that the summands above are 
pure Weil sheaves on $X_0 \times_\kappa \kappa'$, whose $\tau$-weights are all identical.  Suppose each is $\tau$-pure of weight $w$.  For any geometric point $\bar{x}$ over a closed point $x \in X_0$, we set $F_x = F^{[\kappa(x):\kappa]}$; then for any eigenvalue $\alpha$ of $F_x^n$ on an irreducible summand of ${\mathcal G}_{\bar{x}}$, we have $|\tau(\alpha)| = q^{nd(x)w/2}$.  A similar statement thus holds for $F_x$ replacing $F_x^n$ and $q^{wd(x)/2}$ replacing $q^{wd(x)n/2}$.  This holds for every choice of $x$, and $w$ is independent of that choice.  It follows that ${\mathcal G}$ is pure of weight $w$.
\qed

\begin{kor} \label{descend}
Let $r,X$ be as above.  Let ${\mathcal G} \in D^{b,Weil}_c(X_s, \bar{\mathbb Q}_\ell)$.  If $r^*{\mathcal G}$ is a mixed complex in $D^{b,Weil}_c(X_{s'}, \bar{\mathbb Q}_\ell)$, then ${\mathcal G}$ is also mixed.
\end{kor}

\noindent{\em Proof.}
By passing to cohomology sheaves, we can reduce to the case where ${\mathcal G}$ and $r^*{\mathcal G}$ are constructible Weil sheaves on $X_s$ and $X'_{s'}$, respectively.  By Noetherian induction and gluing of constructible sheaves, we can reduce further to the case where ${\mathcal G}$ and thus $r^*{\mathcal G}$ is a smooth Weil sheaf, where the result was proved in the foregoing lemma.
\qed

\subsection{A variant}

The following variant applies to the nearby cycles appearing in the local models and deformations of affine flag varieties treated in this paper.  It does not apply to Shimura varieties: while the generic fibers of these are smooth, they are usually not geometrically integral.

\begin{thm} \label{variant}
Suppose $X_\eta$ is geometrically integral, with intersection complex $^0IC(X_\eta)$.  Then 
$R\Psi^X(\, ^0IC(X_\eta))$ is mixed.
\end{thm}

Here, $^0IC(X_\eta)$ denotes the intersection complex of $X_\eta$, suitably Tate-twisted and shifted to make it a pure self-dual perverse sheaf of weight zero.  This normalization is convenient but of course not necessary for the theorem.

\noindent {\em Proof.}  The proof works the same way as before.  First, replacing $X$ with the scheme-theoretic closure of $X_\eta$, we may assume $X/S$ is flat, hence integral, i.e., we may assume $X$ is an $S$-variety with geometrically integral generic fiber.

Now consider the alterations $\pi: X' \rightarrow X_{S'}$ and $\pi : X'_{\eta'} \rightarrow X_{\eta'}$ associated to the diagram in Theorem \ref{alterations}.  Let $r:X_{S'} \rightarrow X$ denote the projection morphism.  Arguing almost as before (see the Lemma below) we can use the trace map to prove that $^0IC(X_{\eta'}) = \, ^0IC(X_\eta)_{\eta'}$ is a subquotient of $^p\pi_*A^{X'}$.   
By invariance of nearby cycles under change of traits, $r^* R\Psi^{X}(\, ^0IC(X_\eta)) = R\Psi^{X_{S'}}(\, ^0IC(X_{\eta'}))$, which is a subquotient of $R\Psi^{X_{S'}}(\, ^p\pi_* A^{X'}) = \, ^p\pi_*R\Psi^{X'}(A^{X'})$.  The latter is mixed by section 10.2.  Applying descent as in section 10.5, the result follows.
\qed

\begin{lem} \label{subquotient}
The perverse sheaf $\, ^0IC(X_{\eta'})$ is a subquotient of $^p\pi_*(A^{X'})$, in the category $P(X_{\eta'},\Ql)$.
\end{lem}

\noindent {\em Proof.}  Choose the open set $U \subset X_{\eta'}$ as in Lemma \ref{subobject}.  Let $j' : U' := \pi^{-1}(U) \rightarrow X'_{\eta'}$ be the open immersion.  Note that 
$IC := \, ^0IC(X_{\eta'})$ is just $j_{!*} A^{U}$, and $A^{X'} = j'_{!*} A^{U'}$, where $A^{U}$, 
resp. $A^{U'}$ has the obvious meaning.

By using adjunction maps and the morphism of functors $j_{!*} \rightarrow j_*$, we get morphisms 
$$
IC \rightarrow \, ^p\pi_* \pi^* j_{!*}A^U \rightarrow \, ^p\pi_* \pi^* j_*  A^U \rightarrow \, ^p\pi_* j'_* \pi^* A^{U} = \, ^p\pi_* j'_* A^{U'},
$$
and 
$$
^p\pi_*A^{X'} = \, ^p\pi_* j'_{!*}A^{U'} \rightarrow \, ^p \pi_*j'_* A^{U'}.$$
Now restricting to $U$ with $j^*$ and using the trace map as before shows that the first map is non-zero, hence a monomorphism.  Also, restricting shows that the image of $IC$ under the first map lies in the image of the second map (if not, then $IC$ embeds into the cokernel of the second map, which is ruled out after restricting both to $U$).  Thus, $IC$ is a subquotient of $^p\pi_* A^{X'}$.
\qed


\begin{center}
\Large
Part II: A bound on the weights of nearby cycles
\end{center}

The method of proof of Theorem \ref{variant} gives an estimate on the weights:
$$
w(R\Psi^X(\, ^0IC(X_\eta))) \leq {\rm dim}(X_\eta).
$$
(One is reduced to the proper strictly semi-stable case, where it is easy).  

\medskip

Now let $R\Psi = R\Psi(IC(\bar{\mathcal Q}_\mu)[\ell(t_\mu)])$, as usual.

The fact (proved in section 7.2) that for $\mu$ either minuscule, or a coweight for $GL_n$, the function 
${\rm Tr}({\rm Fr}_q, R\Psi_x)$ is polynomial in $q$ of degree $\leq \ell(\mu) - \ell(x)$, is equivalent to the following much stronger bound on the weights of $^0R\Psi := R\Psi(\frac{\ell(\mu)}{2})$: for each $x \in {\rm Adm}(\mu)$, and each closed point $z \in {\mathcal B}_x$, the weights of the complex $^0R\Psi_z$ are $\leq \ell(\mu) - \ell(x)$.  In all the examples we computed, it is true that the trace of Frobenius is always a polynomial in $q$ with this degree.  In all the cases we have proved this (thus far), it follows from the existence of minimal expressions for the functions $\Theta_\lambda$ in the Iwahori-Hecke algebra (cf. section 9.4).  It would be interesting to find a direct geometric argument to prove this sharper bound on the weights in all cases.  This would give a way to deduce the polynomial nature of ${\rm Tr}({\rm Fr}_q, R\Psi_x)$ in general, 
even when the functions $\Theta_\lambda$ are not known to have minimal expressions.

\medskip

The following gives a sufficient condition for the sharp bound to hold.  In this statement, we assume there is the action of a geometrically connected smooth group scheme $G/S$ on $X/S$, and we say a $G_{\bar{s}}$-equivariant Weil perverse sheaf ${\mathcal K}$ on $X_{\bar s}$ satisfies the sharp bound if for each locally closed immersion $i: \mathcal O \rightarrow X_{\bar{s}}$ of a $G_{\bar s}$-orbit ${\mathcal O}$, we have $w(i^*{\mathcal K}) \leq {\rm dim}(X_{\bar{s}}) - {\rm dim}(\mathcal O)$.  

\begin{stz} \label{sharp_bound}
Let $X/S$ be an $S$-variety such that $X_\eta$ is geometrically integral.  Suppose that $X/S$ carries an action of a geometrically connected smooth group scheme $G/S$.  Suppose that there exists a trait $S'$, finite over $S$, and an $S$-alteration $X' \rightarrow X$ with compactification $\overline{X'}$ as in Theorem \ref{alterations}, each carrying an action of $G_{S'}$, such that the inclusion $j: X' \hookrightarrow \overline{X'}$ and the $S'$-alteration  $\pi: X' \rightarrow X_{S'}$ are $G_{S'}$-equivariant.  Then $R\Psi^{X}(\, ^0IC(X_\eta))$ satisfies the sharp bound on its weights.
\end{stz}

\noindent{\em Proof.}  This is similar to the proof of Theorem \ref{variant}.  The sharp bound holds for the strictly semi-stable case, and thus also for $R\Psi^{X'}(A^{X'})$.  It is enough to prove that the sharp bound holds for $^p\pi_!R\Psi^{X'}(A^{X'})$.  So, let ${\mathcal O} \subset X_{s'}$ denote a $G_{s'}$-orbit of dimension $m$, and suppose $z \in {\mathcal O}$ is a closed point.  Using equivariance, one can show that any $G_{s'}$-orbit contained in $\pi^{-1}({\mathcal O})$  dominates ${\mathcal O}$, hence has dimension $\geq m$.  Thus, the weights of 
$\pi_!R\Psi^{X'}(A^{X'})_z = R\Gamma_c(\pi^{-1}(z), R\Psi^{X'}(A^{X'}))$ are $\leq {\rm dim}(X'_{\bar{s}'}) - m$, and we are done.
\qed

There are $\cB$-equivariant semi-stable resolutions for the local models ${\bf M}^{\rm loc} = M_\mu$, in the cases $GSp_4 $, $\mu = (1,1,0,0)$ \cite{deJ1}, and $GSp_6$, $\mu = (1,1,1,0,0,0)$ \cite{Ge}.  Together with the above considerations, we get the following corollary.

\begin{kor}
Let $\mu$ be the minuscule coweight $(1,1,0,0)$ for $GSp_4$ or the minuscule 
coweight $(1,1,1,0,0,0)$ for $GSp_6$.  Then the sharp bound holds for the weights of $^0R\Psi$.  In particular, in each of these cases the multiplicity functions
$$
\sum_i m(R\Psi, w,i) q^i
$$
are polynomials of degree at most $\ell(\mu) - \ell(w)$, for every $w \in {\rm Adm}(\mu)$.
\end{kor}

This Corollary also follows from results already announced in section 9, but the approach given here is perhaps more conceptual.

\bigskip
\bigskip

\noindent {\em Acknowledgements}  This work was initiated during our visits to the Institute for Advanced Study in 2000-2001.  We are grateful to the IAS for its hospitality and financial support.  The research of the second author was also partially supported by NSF grant DMS-0303605.

\bigskip
\bigskip

Ulrich G\"{o}rtz: Mathematisches Institut der Universit\"{a}t Bonn,
Beringstr. 6, 53115 Bonn, Germany.
Email: ugoertz@math.uni-bonn.de

\bigskip

Thomas Haines: Mathematics Department, University of Maryland, College Park, MD 20742-4015, USA. 
Email: tjh@math.umd.edu

\end{document}